# A structure preserving Lanczos algorithm for computing the optical absorption spectrum[*]

Meiyue Shao[1], Felipe H. da Jornada[2,3], Lin Lin[1,4], Chao Yang[1], Jack Deslippe[5], and Steven G. Louie[2,3]

[1]*Computational Research Division, Lawrence Berkeley National Laboratory, Berkeley, CA 94720*
[2]*Department of Physics, University of California, Berkeley, CA 94720*
[3]*Materials Sciences Division, Lawrence Berkeley National Laboratory, Berkeley, CA 94720*
[4]*Department of Mathematics, University of California, Berkeley, CA 94720*
[5]*NERSC, Lawrence Berkeley National Laboratory, Berkeley, CA 94720*

September 5, 2017

## Abstract

We present a new structure preserving Lanczos algorithm for approximating the optical absorption spectrum in the context of solving full Bethe–Salpeter equation without Tamm–Dancoff approximation. The new algorithm is based on a structure preserving Lanczos procedure, which exploits the special block structure of Bethe–Salpeter Hamiltonian matrices. A recently developed technique of generalized averaged Gauss quadrature is incorporated to accelerate the convergence. We also establish the connection between our structure preserving Lanczos procedure with several existing Lanczos procedures developed in different contexts. Numerical examples are presented to demonstrate the effectiveness of our Lanczos algorithm.

---

[*]This work is supported by the Scientific Discovery through Advanced Computing (SciDAC) Program on Excited State Phenomena in Energy Materials funded by the U. S. Department of Energy, Office of Science, under SciDAC program, the Offices of Advanced Scientific Computing Research and Basic Energy Sciences contract number DE-AC02-05CH11231 at the Lawrence Berkeley National Laboratory. It is also partially supported by the Center for Computational Study of Excited-State Phenomena in Energy Materials (C2SEPEM) at the Lawrence Berkeley National Laboratory, which is funded by the U. S. Department of Energy, Office of Science, Basic Energy Sciences, Materials Sciences and Engineering Division under Contract No. DE-AC02-05CH11231, as part of the Computational Materials Sciences Program.

# 1  Introduction

Optical absorption and emission processes provide invaluable information to characterize the electronic properties of solids and molecules. At the same time, an accurate microscopic theory is also highly valuable to predict optical behavior of materials and help design more efficient photovoltaic and light-emitting devices. Physically, the optical spectra of materials can be understood in terms of correlated electron–hole pairs known as excitons. When a photon gets absorbed by a molecule or solid, an electron can be promoted from an occupied to an unoccupied state [25, 31] in a process that creates both a negatively charged particle (known as quasielectron, or simply electron), and a positively charged particle (known as quasihole, or hole). The excitation energy required to produce such an electron–hole pair, or exciton, is directly related to the optical absorption and emission spectrum of the material. A two-particle collective excitation can be described by a two-particle Green's function, of which the real part of the poles give excitation energies. Since the two-particle Green's function satisfies the so called *Bethe–Salpeter equation* (BSE) [26, 31], the excitation energies can be obtained by solving the Bethe–Salpeter equation.

Under an appropriate discretization scheme, the *Bethe–Salpeter Hamiltonian* (BSH) matrix, which is a finite dimensional representation of the Bethe–Salpeter Hamiltonian operator, has the block structure

$$H = \begin{bmatrix} A & B \\ -\overline{B} & -\overline{A} \end{bmatrix} \in \mathbb{C}^{2n \times 2n}, \tag{1}$$

where

$$A^{\mathsf{H}} = A, \qquad \overline{B}^{\mathsf{H}} = B. \tag{2}$$

We can rewrite $H$ as $H = C_n \Omega$ where

$$C_n = \begin{bmatrix} I_n & 0 \\ 0 & -I_n \end{bmatrix}, \qquad \Omega = \begin{bmatrix} A & B \\ \overline{B} & \overline{A} \end{bmatrix}. \tag{3}$$

For most physical systems, $\Omega$ is Hermitian positive definite (see, e.g., [39]), which we will denote by

$$\Omega \succ 0. \tag{4}$$

We define a BSH matrix $H$ that satisfies the condition (4) a *definite* Bethe–Salpeter Hamiltonian matrix. Throughout this paper, we assume that the BSH matrix $H$ is definite, that is, (4) is always assumed.

The matrices $A$ and $B$ are of size $n = n_v n_c n_k$, where $n_v$, $n_c$, and $n_k$ are the numbers of valence states, conduction states, and $k$-points sampled from the Brouillon zone associated with a solid, respectively. Both $n_v$ and $n_c$ are proportional to the number of electrons $n_e$ in the system. Therefore, the dimension $n = O(n_e^2 n_k)$ can be very large for systems of practical interest.

The optical absorption spectrum of a material, which can be measured in spectroscopy experiments, characterizes the excited states properties of the material. Mathematically, the optical absorption spectrum is a matrix functional of the form $d_r^{\mathsf{H}} f(H;\omega) d_l$, where $f(H;\omega)$ is a function of $H$ and the frequency $\omega$, and $d_l$, $d_r \in \mathbb{C}^{2n}$. The function $f(H;\omega)$ has local maxima or minima at eigenvalues of $H$. The locations of these local maxima and minima correspond to excitation and deexcitation energies.

The optical absorption spectrum can be computed by fully diagonalizing the BSH matrix [27]. If $f$ is taken to be the Dirac-$\delta$ distribution, the computed absorption spectrum can be written as

$$\varepsilon_2(\omega) = \sum_{j=1}^{2n} \text{sign}(\lambda_j) \tau_j^2 \delta(\omega - \lambda_j),$$

where $\lambda_j$'s are the eigenvalues of the BSH matrix, and $\tau_j^2$ is known as the oscillator strength associated with $\lambda_j$, and describes the likelihood of excitation (or de-excitation) occurring at energy level $\lambda_j$.

As we will show in the next section, the special block structure of (1) ensures that its eigenvalues appear in positive and negative pairs. Due to the additional structures satisfied by the vectors $d_l$ and $d_r$, the signs of oscillator strengths also appear in positive and negative pairs such that $\varepsilon_2(\omega)$ is of the form

$$\varepsilon_2(\omega) = \sum_{j=1}^{n} \tau_j^2 \big[\delta(\omega - \lambda_j) - \delta(\omega + \lambda_j)\big], \qquad (\tau_j^2 \geq 0).$$

When the problem size $n$ grows, diagonalizing the BSH matrix, whose complexity is $O(n^3)$, becomes increasingly expensive and eventually unaffordable. However, in many cases, the absorption spectrum is simply used as a screening tool for selecting materials with a desirable absorption profile. In these cases, there is no need to accurately locate all eigenvalues and the corresponding oscillator strength—a good approximation to a smoothed absorption spectrum of the form

$$\varepsilon_2(\omega) = \sum_{j=1}^{n} \tau_j^2 \big[\hat{f}(\omega - \lambda_j) - \hat{f}(\omega + \lambda_j)\big], \qquad (\tau_j^2 \geq 0),$$

where $\hat{f}$ is a smooth approximation to the Dirac-$\delta$ distribution, is often sufficient. The use of smooth approximation to the Dirac-$\delta$ distribution is physically meaningful because electron excitation has a finite life-time, and the width of an isolated peak in $\hat{f}$ is related to the inverse of the life time. In this paper, we discuss how to estimate the absorption spectrum efficiently and reliably without diagonalizing $H$.

Our basic idea is to use a $k$-step iterative method to construct a Krylov subspace $\mathcal{S}$ onto which the BSH is projected with $2k \ll n$. We compute the eigenvalues $\theta_j$, $j = 1$, 2, ..., 2k of the $2k \times 2k$ projected Hamiltonian (i.e., the Ritz values), and use them to construct an approximate absorption spectrum of the form

$$\varepsilon_2(\omega) \approx \sum_{j=1}^{k} \tilde{\tau}_j^2 \big[\hat{f}(\omega - \theta_j) - \hat{f}(\omega - \theta_{j+k})\big] = \sum_{j=1}^{k} \tilde{\tau}_j^2 \big[\hat{f}(\omega - \theta_j) - \hat{f}(\omega + \theta_j)\big], \quad (5)$$

where $\tilde{\tau}_j^2 \geq 0$ are approximate oscillator strengths extracted from the same subspace. In order to obtain an approximation of the form (5), we need to exploit properties of $H$ to preserve key properties of (5) when constructing the subspace $\mathcal{S}$. In particular,

- we would like the Ritz values to be real and appear in positive and negative pairs;

- we would like to ensure that the signs of approximate oscillator strengths also appear in positive and negative pairs, and are nonnegative for positive Ritz values.

An algorithm that achieves the above two criteria when it is used to construct $\mathcal{S}$ is called a *structure preserving* algorithm.

In the context of *Tamm–Dancoff approximation* (TDA) [8, 32], which sets the off-diagonal blocks $B$ in $H$ to zero, we only need to consider the matrix $A$ in $H$, which is Hermitian. As a result, we can use the Lanczos algorithm to construct a Krylov subspace from which a structure preserving approximate absorption spectrum can be easily obtained. We can choose either a Gaussian function or a Lorentzian function as $\hat{f}$. The latter choice is adopted in the so called Haydock's algorithm [14, 37].

For full BSE calculations, which is non-Hermitian, the Arnoldi or nonsymmetric Lanczos algorithms [2, 15, 21, 34, 35] do not produce a structure preserving subspace or a projected Hamiltonian that satisfies the desirable properties listed earlier.

Recently a special Lanczos algorithm applicable to a full BSE has been proposed in [13]. In this paper, we analyze additional properties of this Lanczos algorithm and develop an alternative algorithm that works at least equally well. The motivation is to show that both algorithms can produce approximate absorption spectra of the form (5), and are hence structure preserving. One main advantage of the new algorithm presented in this paper is that it allows us to incorporate a new generalized Gauss quadrature rule [19, 30] to further improve the accuracy of the approximate absorption spectrum with negligible additional cost. We examine the connection between the new algorithm and the algorithm presented in [13], as well as a few other variants of the Lanczos algorithms.

The rest of the paper is organized as follows. In Section 2, we review some basic properties of the definite Bethe–Salpeter Hamiltonian and the optical absorption spectrum. In Section 3, we describe how the standard Lanczos algorithm can be used to estimate the absorption spectrum in the context of TDA. Then, in Section 4, we discuss how the Lanczos algorithm can be modified to preserve the desirable structures of BSH and the corresponding absorption spectrum to be approximated. We compare several variants of the Lanczos algorithm, examine the connection among them, and discuss whether they preserve the desirable properties of the absorption spectrum to be approximated. Finally, computational examples are presented in Section 5 to demonstrate the effectiveness and efficiency of the proposed Lanczos algorithm.

## 2 Preliminaries

### 2.1 Properties of definite Bethe–Salpeter Hamiltonian matrices

We first briefly review some basic spectral properties of definite BSH matrices. Detailed discussion on these properties can be found in [5, 27, 29].

Although a definite BSH matrix $H$ defined in (1) is in general non-Hermitian, it is diagonalizable and has real spectrum. Moreover the special structure of the BSH leads to a structured spectral decomposition as stated in Theorem 1 below.

**Theorem 1** ([27, Theorem 3]). *Let $H$ be a definite Bethe–Salpeter Hamiltonian matrix. Then the spectral decomposition of $H$ is of the form $H = Z \operatorname{diag}\{\Lambda_+, \Lambda_-\} Z^{-1}$ where*

$$Z = \begin{bmatrix} X & \overline{Y} \\ Y & \overline{X} \end{bmatrix}, \qquad Z^{-1} = C_n Z^{\mathsf{H}} C_n = \begin{bmatrix} X & -\overline{Y} \\ -Y & \overline{X} \end{bmatrix}^{\mathsf{H}}, \tag{6}$$

$\Lambda_+ = \operatorname{diag}\{\lambda_1, \lambda_2, \ldots, \lambda_n\}$, *and* $\Lambda_- = \operatorname{diag}\{\lambda_{n+1}, \lambda_{n+2}, \ldots, \lambda_{2n}\}$ *with*

$$\lambda_1 = -\lambda_{n+1} \geq \lambda_2 = -\lambda_{n+2} \geq \cdots \geq \lambda_n = -\lambda_{2n} > 0.$$

Since the eigenvalues of $H$ appear in positive and negative pairs $\pm\lambda_j$, we use $\lambda_j^+ \equiv \lambda_j$ and $\lambda_j^- \equiv -\lambda_j$, for $1 \leq j \leq n$, in the following to emphasize on the signs of these eigenvalues. Let $X = [x_1, \ldots, x_n]$, $Y = [y_1, \ldots, y_n] \in \mathbb{C}^{n \times n}$ be the submatrices in (6). Theorem 1 suggests that the right and left eigenvectors associated with the positive eigenvalue $\lambda_j^+$ are $z_j = [x_j^{\mathsf{H}}, y_j^{\mathsf{H}}]^{\mathsf{H}}$ and $C_n z_j = [x_j^{\mathsf{H}}, -y_j^{\mathsf{H}}]^{\mathsf{H}}$ respectively, and the right and left eigenvectors associated with $\lambda_j^-$ are $z_{n+j} = [\overline{y}_j^{\mathsf{H}}, \overline{x}_j^{\mathsf{H}}]^{\mathsf{H}}$ and $-C_n z_{n+j} = [-\overline{y}_j^{\mathsf{H}}, \overline{x}_j^{\mathsf{H}}]^{\mathsf{H}}$ respectively. The normalization condition $(C_n Z C_n)^{\mathsf{H}} Z = I_{2n}$ implies that

$$x_j^{\mathsf{H}} x_j - y_j^{\mathsf{H}} y_j = 1$$

for $j = 1, \ldots, n$. As long as the right eigenvectors associated with the positive eigenvalues are properly normalized, other eigenvectors can be easily recovered.

From (6), we observe that the right eigenvectors of $H$ are orthonormal with respect to the $C$-inner product, $\langle u, v \rangle_{C_n} = v^{\mathsf{H}} C_n u$, which is an *indefinite inner product*. Another observation is

$$Z^{\mathsf{H}} \Omega Z = Z^{\mathsf{H}} C_n Z \operatorname{diag}\{\Lambda_+, -\Lambda_+\} = C_n \operatorname{diag}\{\Lambda_+, -\Lambda_+\} = \operatorname{diag}\{\Lambda_+, \Lambda_+\}, \tag{7}$$

indicating that the right eigenvectors of $H$ are also orthogonal with respect to the $\Omega$-inner product $\langle u, v \rangle_\Omega = v^{\mathsf{H}} \Omega u$. These orthogonalities are crucial for developing structure preserving Lanczos procedures. By structure preserving, we mean that the positive and negative pairing of the eigenvalues is preserved in the approximations to the eigenvalues of BSH.

## 2.2 Optical absorption spectra

Let $(z_r)_j$ and $(z_l)_j$ be the right and left eigenvectors of $H$, respectively, associated with the eigenvalue $\lambda_j$, $(1 \leq j \leq 2n)$. We denote by $\varepsilon_2(\omega)$ the imaginary part of the macroscopic dielectric function; $\varepsilon_2(\omega)$ is also proportional to the optical absorption spectrum of a material, and can be computed in a straightforward way from the eigenvalues and eigenvectors of the BSH as

$$\varepsilon_2(\omega) = \frac{8\pi^2 e^2}{V_{\text{xtal}}} \epsilon(\omega),$$

$$\epsilon(\omega) := d_r^{\mathsf{H}} \delta(\omega I_{2n} - H) d_l = \sum_{j=1}^{2n} \frac{(d_r^{\mathsf{H}}(z_r)_j)((z_l)_j^{\mathsf{H}} d_l)}{(z_l)_j^{\mathsf{H}}(z_r)_j} \delta(\omega - \lambda_j), \quad (8)$$

where $V_{\text{xtal}}$ is the crystal volume, $e$ is the elementary charge, and

$$d_r = \begin{bmatrix} d \\ -\overline{d} \end{bmatrix} \quad \text{and} \quad d_l = \begin{bmatrix} d \\ \overline{d} \end{bmatrix}$$

are the right and left *optical transition vectors*, respectively. Because the $d_r$ and $d_l$ depend solely on $d$, we will simply refer to $d$ as the optical transition vector. The coefficient, $(d_r^{\mathsf{H}}(z_r)_j)((z_l)_j^{\mathsf{H}} d_l)/((z_l)_j^{\mathsf{H}}(z_r)_j)$, of the Dirac delta function $\delta(\omega - \lambda_j)$ in (8) is known as the *oscillator strength* associated to the excitonic state $j$. Figure 1 shows a typical curve for the imaginary part of the dielectric function. In order to produce this plot, each Dirac delta function was broadened by a Gaussian function, as we will discuss below. The height of each peak in the spectrum is determined by the oscillator strength associated to each eigenvalue $\lambda_j$ and the number of eigenvalues clustered around an energy. Since the optical absorption spectrum is proportional to $\varepsilon_2(\omega)$, which is in turn proportional to $\epsilon(\omega)$, in this work we will broadly refer to both $\epsilon(\omega)$ and $\varepsilon_2(\omega)$ as the optical absorption spectrum of a material.

If $H$ can be fully diagonalized, we can compute $\epsilon(\omega)$ using the eigenpairs of $H$. However, diagonalizing $H$ is often costly, especially when the dimension of $H$ becomes large. For instance, for many low dimensional systems such as monolayer $MoS_2$, $n$ is on the order of 360,000 [18]. Similarly, large $n$'s are required to fully converge calculations on semiconducting carbon nanotubes and to obtain the correct order of the excited excitonic states in bulk semiconductors, for example. Therefore, it is natural to seek alternative approaches.

Using the structure of the eigenvectors of the BSH matrix $H$, we can simplify the expression of the absorption spectrum. For positive eigenpairs, we can choose $(z_r)_j = z_j$

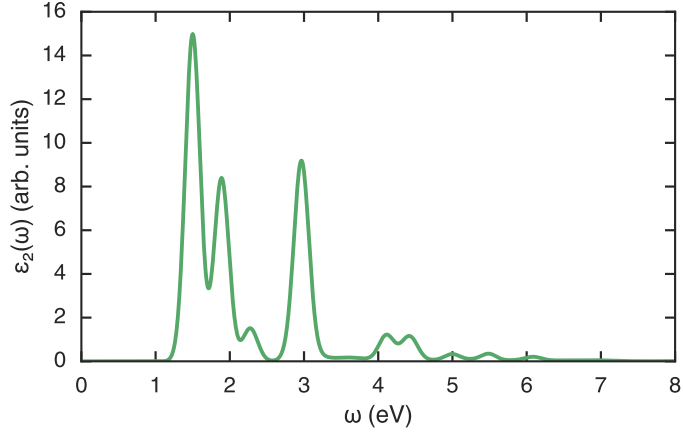

Figure 1: A typical curve for the imaginary part of the dielectric function. This curve is obtained from a single-wall $(8,0)$ carbon nanotube.

and $(z_l)_j = C_n z_j$ so that $(z_l)_j^{\mathsf{H}}(z_r)_j = 1$. Then we have

$$(z_l)_j^{\mathsf{H}} d_l = (C_n z_j)^{\mathsf{H}}(C_n d_r) = z_j^{\mathsf{H}} d_r = \overline{d_r^{\mathsf{H}} z_j}.$$

It follows that

$$\epsilon^+(\omega) := \sum_{j=1}^{n} \frac{(d_r^{\mathsf{H}}(z_r)_j)((z_l)_j^{\mathsf{H}} d_l)}{(z_l)_j^{\mathsf{H}}(z_r)_j} \delta(\omega - \lambda_j)$$

$$= \sum_{j=1}^{n} |d_r^{\mathsf{H}} z_j|^2 \delta(\omega - \lambda_j^+)$$

$$= \sum_{j=1}^{n} |d^{\mathsf{H}} x_j - \overline{d}^{\mathsf{H}} y_j|^2 \delta(\omega - \lambda_j^+).$$

We remark that the oscillator strength $|d_r^{\mathsf{H}} z_j|^2$ is nonnegative. Similarly, for negative eigenpairs, we have

$$(z_l)_j^{\mathsf{H}} d_l = (-C_n z_j)^{\mathsf{H}}(C_n d_r) = -z_j^{\mathsf{H}} d_r = -\overline{d_r^{\mathsf{H}} z_j}$$

and

$$\epsilon^-(\omega) := -\sum_{j=n+1}^{2n} |d_r^{\mathsf{H}} z_j|^2 \delta(\omega - \lambda_j) = -\sum_{j=1}^{n} |d^{\mathsf{H}} x_j - \overline{d}^{\mathsf{H}} y_j|^2 \delta(\omega + \lambda_j^+) = -\epsilon^+(-\omega).$$

Therefore, the absorption spectrum

$$\epsilon(\omega) = \epsilon^+(\omega) + \epsilon^-(\omega) = \epsilon^+(\omega) - \epsilon^+(-\omega)$$

can be viewed as an odd function of the frequency $\omega$ in the distribution sense.

In practice, it is not desirable to plot the imaginary part of the polarizability as a sum of Dirac delta functions. A broadened peaked function, such as the Lorentzian function

$$L_\sigma(\omega) := \frac{1}{\pi} \cdot \frac{\sigma}{\omega^2 + \sigma^2} = \frac{1}{\pi} \mathrm{Im} \frac{1}{\omega - i\sigma}$$

or the Gaussian function

$$G_\sigma(\omega) := \frac{1}{\sqrt{2\pi}\,\sigma} e^{-\omega^2/(2\sigma^2)},$$

is used to replace the Dirac delta function, where the broadening factor $\sigma > 0$ is a small number. The first reason for doing so is because there is physically a lifetime associated to each excitonic state. The second reason is due to discretization procedures in performed in the calculations, such as employing a finite number of $k$-points in calculations on extended systems. If a calculation could be carried with infinitely many $k$-points, the optical absorption spectrum would consist of a few isolated low-energy sharp peaks, but the delta functions merge at higher energy and form a continuum spectrum. On the other hand, a calculation performed with a finite number of $k$-points only samples a finite number of transitions in this continuum region.

Therefore, we wish to plot the imaginary part of the dielectric function using a generic peaked function $g_\sigma(\omega)$ (either $L_\sigma(\omega)$ or $G_\sigma(\omega)$) characterized by a typical small width $\sigma$ instead of $\delta(\omega)$. The imaginary part of the dielectric matrix can be expressed now in terms of

$$\epsilon_\sigma(\omega) = d_r^{\mathsf{H}} g_\sigma(\omega I_{2n} - H) d_l = \sum_{j=1}^n \left| d^{\mathsf{H}} x_j - \overline{d}^{\mathsf{H}} y_j \right|^2 \left[ g_\sigma(\omega - \lambda_j^+) - g_\sigma(\omega + \lambda_j^+) \right], \qquad (9)$$

which is an odd function in $\omega$, that is, $\epsilon_\sigma(-\omega) = -\epsilon_\sigma(\omega)$. Thus it suffices to compute the function value for $\omega > 0$.

Note that (9), which is a scalar function of $\omega$, can be viewed as an expected value of a matrix function. We are interested in the positions and heights of the major peaks of this function, which are given by the eigenvalues and eigenvectors of the BSH. However, the precise position and height of each peak is seldom required, especially since the underlying theories employed to obtain these spectra are already themselves approximate. Therefore, efficient methods that can provide estimates of (9) without computing each individual eigenpair of $H$ are of great interest. In Sections 3 and 4, we discuss how to use Lanczos algorithms to estimate $\epsilon_\sigma(\omega)$ efficiently.

## 3   Tamm–Dancoff approximation

Tamm–Dancoff approximation (TDA) [8, 25, 32] is a technique often used to in practice to reduce the computational cost of the absorption spectrum calculation. For many systems, especially on bulk semiconductors and metals, the TDA incurs a very small error in the optical absorption spectrum, and for that reason it has been a widely used approximation in condensed-matter physics. In this section we discuss how to estimate the absorption spectrum with a Lanczos procedure within the TDA.

However, we remark that for many systems, including systems with reduced dimensionality optically excited with light polarized along a confined direction, the TDA may incur in large errors for the optical absorption spectrum. We shall discuss full BSE solvers in Section 4.

### 3.1 Lanczos algorithm

In TDA, the off-diagonal block of $H$, $B$, is set to zero. We denote the resulting block diagonal BSH by $H_{\text{TDA}} = \text{diag}\{A, -\overline{A}\}$, which is a Hermitian matrix. It follows that the absorption spectrum associated with $H_{\text{TDA}}$ becomes

$$\epsilon(\omega) = d_r^{\mathsf{H}} \delta(\omega I_{2n} - H_{\text{TDA}}) d_l = d^{\mathsf{H}} \delta(\omega I_n - A) d - \overline{d^{\mathsf{H}} \delta(\omega I_n + A) d}.$$

As $d^{\mathsf{H}} \delta(\omega I_n \pm A) d$ is real and nonnegative, we can omit the complex conjugation in the second term. In practice, we compute

$$\epsilon_\sigma(\omega) = d^{\mathsf{H}} g_\sigma(\omega I_n - A) d - d^{\mathsf{H}} g_\sigma(\omega I_n + A) d =: d^{\mathsf{H}} f(A; \omega) d \tag{10}$$

for $\omega > 0$, where $f(t; \omega) = g_\sigma(\omega - t) - g_\sigma(\omega + t)$.

Since $A$ is Hermitian and positive definite, the matrix functionals in (10) can be estimated using the Lanczos algorithm. Starting with $u_1 = d/\|d\|_2$, a $k$-step Lanczos procedure produces

$$AU_k = U_k T_k + \beta_k u_{k+1} e_k^{\mathsf{H}}, \qquad U_{k+1}^{\mathsf{H}} U_{k+1} = I_{k+1}, \tag{11}$$

where

$$T_k = \text{tridiag} \left\{ \begin{matrix} & \beta_1 & \cdots & & \beta_{k-1} & \\ \alpha_1 & \cdots & & \cdots & & \alpha_k \\ & \beta_1 & \cdots & & \beta_{k-1} & \end{matrix} \right\} \tag{12}$$

is a real symmetric, tridiagonal, positive definite, and componentwise nonnegative matrix. Here we use the convention $U_j = [u_1, \ldots, u_j]$ (for $1 \leq j \leq k+1$) to represent the Lanczos vectors, and $e_j$ is the $j$th column of the identity matrix. Then the absorption spectrum can be estimated by

$$d^{\mathsf{H}} f(A; \omega) d = \|d\|_2^2 u_1^{\mathsf{H}} f(A; \omega) U_k e_1 \approx \|d\|_2^2 u_1^{\mathsf{H}} U_k f(T_k; \omega) e_1 = \|d\|_2^2 e_1^{\mathsf{H}} f(T_k; \omega) e_1. \tag{13}$$

As long as $k \ll n$, the matrix function of the projected matrix $T_k$, $f(T_k; \omega)$, can be easily evaluated by diagonalizing $T_k$. Moreover, there is no need to explicitly store the whole history of the Lanczos vectors because eventually only $T_k$ is used in (13). However, it is important to ensure columns of the generated $U_{k+1}$ matrix are orthonormal. A desired feature here is that the estimated absorption spectrum in this approach is nonnegative for $\omega > 0$. Clearly, the Lanczos algorithm possesses this desired feature.

Finally, we remark that when the Gaussian functions are replaced by Lorentzian functions Haydock's recursive algorithm is mathematically equivalent to the Lanczos algorithm. As the Lanczos algorithm is more general—it can handle any approximation to the Dirac delta function, it is a simple and flexible replacement of Haydock's recursive method in this context. Another advantage of the Lanczos algorithm will be discussed in the next subsection.

## 3.2 Generalized averaged Gauss quadrature

It is well known that the Lanczos algorithm for estimating matrix functionals can be interpreted as Gauss quadrature [10, 11]. In [19], a recently developed generalized averaged Gauss quadrature rule [30] has been adopted to improve the accuracy of the Lanczos algorithm with little extra effort. In the following we briefly describe the procedure of this approach.

After the $k$-step Lanczos procedure is performed, we can construct a $(2k-1) \times (2k-1)$ symmetric tridiagonal matrix $\widehat{T}_k$ as

$$\widehat{T}_k = \text{tridiag} \left\{ \begin{array}{ccccccccc} & \beta_1 & \cdots & \beta_{k-1} & \beta_k & \beta_{k-2} & \cdots & \beta_1 & \\ \alpha_1 & \cdots & \cdots & & \alpha_k & \alpha_{k-1} & \cdots & \cdots & \alpha_1 \\ & \beta_1 & \cdots & \beta_{k-1} & \beta_k & \beta_{k-2} & \cdots & \beta_1 & \end{array} \right\}. \tag{14}$$

Then we replace $e_1^{\mathsf{H}} f(T_k;\omega)e_1$ in (13) by $e_1^{\mathsf{H}} f(\widehat{T}_k;\omega)e_1$,[1] that is,

$$d^{\mathsf{H}} f(A;\omega)d \approx \|d\|_2^2 e_1^{\mathsf{H}} f(\widehat{T}_k;\omega)e_1. \tag{15}$$

When $k$ is not very large, the cost of computing $f(T_k;\omega)$ or $f(\widehat{T}_k;\omega)$ is negligible compared to that of forming $T_k$. As the spectrum of $\widehat{T}_k$ is a superset of that of $T_{k-1}$, and $\Lambda(\widehat{T}_k)\backslash\Lambda(T_{k-1})$ interlaces with $\Lambda(T_{k-1})$, (15) should be a better approximation compared to (13) with negligible computational overhead. We refer the readers to [19, 30] for detailed discussions.

If the Lanczos procedure breaks down at $k$th step, that is, $\beta_k = 0$, then (13) holds exactly instead of approximately. In this lucky breakdown, (15) also holds exactly because $\widehat{T}_k$ decouples into two tridiagonal submatrices. We remark that an extra benefit of using the generalized averaged Gauss quadrature is that, for the same number of quadrature points, the generalized averaged Gauss quadrature requires fewer Lanczos steps, hence the risk of loss of orthogonality among the Lanczos vectors is reduced.

Certainly the generalized averaged Gauss quadrature can be adopted here for the estimation of absorption spectrum. Similar to the Lanczos algorithm with standard Gauss quadrature, the generalized averaged Gauss quadrature also produces nonnegative oscillator strengths. Thus the estimated absorption spectrum is also nonnegative for $\omega > 0$ when $\widehat{T}_k$ is positive definite. However, $\widehat{T}_k$ as defined in (14) can sometimes have one nonpositive eigenvalue. (The second smallest eigenvalue of $\widehat{T}_k$ is always positive since $\Lambda(\widehat{T}_k)\backslash\Lambda(T_{k-1})$ interlaces with $\Lambda(T_{k-1})$.) Such a nonpositive eigenvalue may violate the property $\epsilon_\sigma(\omega) \geq 0$ for $\omega > 0$. A simple remedy is to redefine $f(t; \omega)$ as

$$f(t; \omega) = \begin{cases} g_\sigma(\omega - t) - g_\sigma(\omega + t) & \text{if } t > 0, \\ 0 & \text{if } t \leq 0. \end{cases}$$

---

[1]The vector $e_1$ is of length $k$ in $e_1^{\mathsf{H}} f(T_k;\omega)e_1$, and is of length $2k-1$ in $e_1^{\mathsf{H}} f(\widehat{T}_k;\omega)e_1$.

**Algorithm 1** Lanczos algorithm for estimating the absorption spectrum under TDA.
***
**Input:** A Hermitian positive definite matrix $A \in \mathbb{C}^{n \times n}$, an optical transition vector $d \in \mathbb{C}^n$, a broadening factor $\sigma > 0$, the number of Lanczos steps $k$, and a set of frequencies $\{\omega_i\}_{i=1}^m$.
**Output:** The estimated absorption spectrum $\epsilon_\sigma(\omega)$ sampled at $\omega_i$ (for $1 \leq i \leq m$).
1: Perform $k$ Lanczos steps using $d$ as the starting vector.
2: Formulate $\widehat{T}_k$ as defined in (14).
3: Compute the spectral decomposition $\widehat{T}_k = \widehat{S}_k \, \mathrm{diag}\left\{\widehat{\theta}_1, \ldots, \widehat{\theta}_{2k-1}\right\} \widehat{S}_k^{\mathsf{H}}$, where $\widehat{S}_k^{\mathsf{H}} \widehat{S}_k = I_{2k-1}$.
4: Evaluate
$$\epsilon_\sigma(\omega_i) = \|d\|_2^2 \sum_{\substack{j=1 \\ \widehat{\theta}_j > 0}}^{2k-1} |\widehat{S}_k(1,j)|^2 \bigl[g_\sigma(\omega_i - \widehat{\theta}_j) - g_\sigma(\omega_i + \widehat{\theta}_j)\bigr]$$
for $i = 1, \ldots, m$.
***

Then in the resulting generalized averaged Gauss quadrature (15) we can simply discard the term involving the nonpositive eigenvalue of $\widehat{T}_k$, if there is any. In fact, dropping the nonpositive eigenvalue does not affect the accuracy, because the eigenvalues of $T_{k-1}$, as the common Gauss quadrature nodes for both (13) and (15) (assuming in (13) we use the approximation from a $(k-1)$-step Lanczos procedure instead of a $k$-step one), have the same weights (up to scaling) in both quadrature rules [30]. We summarize the Lanczos algorithm with generalized averaged Gauss quadrature in Algorithm 1. The utility of generalized averaged Gauss quadrature provides another advantage of the Lanczos algorithm over Haydock's recursive algorithm.

## 4 Absorption spectrum for full BSE

In this section we investigate how to estimate the absorption spectrum without using the Tamm–Dancoff approximation. Like the Lanczos algorithm in the TDA setting, the following features are desired.

1. Any breakdown in the Lanczos procedure is a lucky breakdown.

2. The computed absorption spectrum is real and nonnegative for $\omega > 0$.

3. The full history of Lanczos vectors is not required.

4. The technique of generalized averaged Gauss quadrature can be applied.

We shall demonstrate that all these features are feasible for full BSE calculations.

## 4.1 Lanczos algorithm for real BSE

We first examine a simpler case in which both $A$ and $B$ are real symmetric matrices, and

$$H = \begin{bmatrix} A & B \\ -B & -A \end{bmatrix} \in \mathbb{R}^{2n \times 2n} \tag{16}$$

is real also. Such an $H$ results from systems with real-space inversion symmetry. It is not difficult to verify that the condition (4) is equivalent to the following conditions:

$$M := A + B \succ 0, \qquad K := A - B \succ 0. \tag{17}$$

We also assume that the optical transition vector $d$ is real. A Lanczos algorithm that can be used to estimate $\epsilon(\omega)$ for BSH matrices of this type has been studied in [7], in the context of linear response time-dependent density functional theory (TDDFT) based calculations. In the following, we briefly summarize this algorithm.

Using the spectral decomposition of $H$ as shown in Theorem 1, we can verify that

$$M = (X - Y)\Lambda_+(X - Y)^\mathsf{T}, \qquad K = (X + Y)\Lambda_+(X + Y)^\mathsf{T},$$

and

$$(X - Y)^\mathsf{T}(X + Y) = I_n.$$

Then we have

$$\begin{aligned}
\epsilon(\omega) &= \sum_{j=1}^{n} \left[d^\mathsf{T}(x_j - y_j)\right]^2 \left[\delta(\omega - \lambda_j^+) - \delta(\omega + \lambda_j^+)\right] \\
&= 2\operatorname{sign}(\omega) \sum_{j=1}^{n} \lambda_j^+ \left[d^\mathsf{T}(x_j - y_j)\right]^2 \delta\!\left(\omega^2 - (\lambda_j^+)^2\right) \\
&= 2\operatorname{sign}(\omega) d^\mathsf{T}(X - Y)\Lambda_+(X - Y)^\mathsf{T}(X + Y)\delta(\omega^2 I_n - \Lambda_+^2)(X - Y)^\mathsf{T} d \\
&= 2\operatorname{sign}(\omega) d^\mathsf{T} M \delta(\omega^2 I_n - KM) d.
\end{aligned}$$

**Algorithm 2** Lanczos procedure in $M$-inner product for real full BSE.

**Input:** A definite Bethe–Salpeter Hamiltonian matrix $H \in \mathbb{R}^{2n \times 2n}$; a starting vector $u_1 \in \mathbb{R}^n$ satisfying $u_1^\mathsf{T}(A+B)u_1 = 1$; the number of Lanczos steps, $k$.
**Output:** $\alpha_1, \ldots, \alpha_k, \beta_1, \ldots, \beta_k \in \mathbb{R}$, and $u_1, \ldots, u_{k+1} \in \mathbb{R}^n$ satisfying (18) and $U_{k+1}^\mathsf{T}(A+B)U_{k+1} = I_{k+1}$.
 1: $\beta_0 \leftarrow 0, \quad u_0 \leftarrow 0, \quad v_0 \leftarrow 0$.
 2: $v_1 = (A+B)u_1$.
 3: **for** $j = 1, \ldots, k$ **do**
 4: $\quad x \leftarrow (A-B)v_j - \beta_{j-1}u_{j-1}$.
 5: $\quad \alpha_j \leftarrow v_j^\mathsf{T} x$.
 6: $\quad x \leftarrow x - \alpha_j u_j$.
 7: $\quad y \leftarrow (A+B)x$.
 8: $\quad \beta_j \leftarrow \sqrt{x^\mathsf{T} y}$.
 9: $\quad u_{j+1} \leftarrow x/\beta_j, \quad v_{j+1} \leftarrow y/\beta_j$.
10: **end for**

Therefore, we reduce this problem size from $2n \times 2n$ to $n \times n$. Although $KM$ is nonsymmetric in general, it is symmetric and positive definite with respect to the $M$-inner product because
$$\langle x, KMy \rangle_M = y^\mathsf{T} MKMx = \langle KMx, y \rangle_M.$$
A Lanczos procedure in which standard Euclidean inner product is replaced with an $M$-inner product reads
$$KMU_k = U_k T_k + \beta_k u_{k+1} e_k^\mathsf{T}, \qquad (18)$$
with $u_1 = d/\|d\|_M$ and $U_{k+1}^\mathsf{T} MU_{k+1} = I_{k+1}$. Algorithm 2 outlines the computational procedure of calculating (18). We remark that in [7] full orthogonalization is used to maintain numerical stability of the Lanczos procedure. In contrast, Algorithm 2 uses a careful formulation of short recurrence. The numerical stability is observed to be comparable with full orthogonalization if $k$ is reasonably small.

It follows from (18) and the identity
$$\delta(\omega - |\lambda|) - \delta(\omega + |\lambda|) = 2|\lambda|\operatorname{sign}(\omega)\delta(\omega^2 - \lambda^2)$$
that $\epsilon(\omega)$ can be approximated through
$$\begin{aligned}
\epsilon(\omega) &= 2\operatorname{sign}(\omega) d^\mathsf{T} M \delta(\omega^2 I_n - KM) d \\
&\approx 2\operatorname{sign}(\omega) d^\mathsf{T} M U_k \delta(\omega^2 I_k - T_k) U_k^\mathsf{T} M d \\
&= \|d\|_M^2 e_1^\mathsf{T} \left[\delta(\omega I_k - T_k^{1/2}) - \delta(\omega I_k + T_k^{1/2})\right] T_k^{-1/2} e_1 \\
&\approx \|d\|_M^2 e_1^\mathsf{T} \left[g_\sigma(\omega I_k - T_k^{1/2}) - g_\sigma(\omega I_k + T_k^{1/2})\right] T_k^{-1/2} e_1.
\end{aligned} \qquad (19)$$

Here $T_k$ is a real symmetric tridiagonal matrix as in (12). Similar to the Lanczos algorithm in the TDA setting, the approximate $\epsilon_\sigma(\omega)$ is nonnegative for $\omega > 0$, which is a desired property. There is also no need to keep the whole history of Lanczos vectors.

We have already seen in Section 3.2 that the generalized averaged Gauss quadrature can be incorporated in the Lanczos algorithm. This is also the case for (19). Let

$$f(t;\omega) = t^{-1/2}\big[g(\omega - t^{1/2}) - g(\omega + t^{1/2})\big].$$

Then the generalized averaged Gauss quadrature replaces $e_1^\mathsf{T} f(T_k;\omega)e_1$ in (19) by $e_1^\mathsf{T} f(\widehat{T}_k;\omega)e_1$, that is,

$$\epsilon_\sigma(\omega) \approx \|d_1\|_M^2 e_1^\mathsf{T} f(\widehat{T}_k;\omega)e_1, \tag{20}$$

where $\widehat{T}_k \in \mathbb{R}^{(2k-1)\times(2k-1)}$ is as defined in (14). It is expected that (20) in general provides a better approximation to $\epsilon_\sigma(\omega)$ compared to (19). Similar to the discussions in Section 3.2, $\widehat{T}_k$ can sometimes has one nonpositive eigenvalue. But in (19) $\widehat{T}_k$ needs to be positive definite so that $\widehat{T}_k^{1/2}$ is also positive definite. The remedy is to extend the definition of $f(t;\omega)$ as

$$f(t;\omega) = \begin{cases} t^{-1/2}\big[g_\sigma(\omega - t^{1/2}) - g_\sigma(\omega + t^{1/2})\big] & \text{if } t > 0, \\ 0 & \text{if } t \leq 0, \end{cases}$$

and discard the term involving the nonpositive eigenvalue of $\widehat{T}_k$, if there is any. Algorithm 3 summarizes the Lanczos algorithm for real full BSE incorporated with the generalized averaged Gauss quadrature.

We should point out that $\epsilon(\omega)$ can be obtained by computing the eigenpairs of $KM$ or $H$ directly. If one is only interested in the low energy region of the absorption spectrum, iterative methods such as the ones proposed in [1, 20, 36] can be used to compute the first few eigenpairs. However, these methods can become costly when the absorption spectrum window becomes large, and more eigenpairs need to be computed.

## 4.2 Structure preserving Lanczos procedure for complex BSE

In this subsection we discuss how to develop a structure preserving Lanczos procedure for complex BSE. Just like the real case, we will try to reformulate the problem so that only $n$-dimensional matrices and vectors are involved. To this end, let us define

$$\mathcal{U}_\phi = \left\{ \begin{bmatrix} u \\ \mathrm{e}^{\mathrm{i}\phi}\overline{u} \end{bmatrix} : u \in \mathbb{C}^n \right\}, \qquad (\phi \in \mathbb{R}).$$

It can be easily verified that $H\mathcal{U}_\phi = \mathcal{U}_{\phi+\pi}$ and $H^2 \mathcal{U}_\phi = \mathcal{U}_\phi$. However, we remark that $\mathcal{U}_\phi$ is *not* an invariant subspaces of $H^2$ as it is not a subspace of $\mathbb{C}^{2n}$ over $\mathbb{C}$; it can only be regarded as a linear space over $\mathbb{R}$. To approximate $d_r^\mathsf{H} \delta(\omega I_{2n} - H)d_l$ using a Lanczos procedure, it is natural to use $d_l$ as the starting vector. Note that $d_l$ and $d_r$ are structured

---
**Algorithm 3** The Lanczos algorithm for estimating the optical absorption spectrum for real full BSE.
---
**Input:** Real symmetric positive definite matrices $M$, $K \in \mathbb{R}^{n \times n}$, an optical transition vector $d \in \mathbb{R}^n$, a broadening factor $\sigma > 0$, the number of Lanczos steps $k$, and a set of frequencies $\{\omega_i\}_{i=1}^m$.
**Output:** The estimated absorption spectrum $\epsilon_\sigma(\omega)$ sampled at $\omega_i$ (for $1 \leq i \leq m$).
1: Perform $k$ Lanczos steps in $M$-inner product using $d$ as the starting vector.
2: Formulate $\widehat{T}_k$ as defined in (14).
3: Compute the spectral decomposition $\widehat{T}_k = \widehat{S}_k \operatorname{diag}\left\{\widehat{\theta}_1^2, \ldots, \widehat{\theta}_{2k-1}^2\right\} \widehat{S}_k^{\mathsf{H}}$, where $\widehat{S}_k^{\mathsf{H}} \widehat{S}_k = I_{2k-1}$ and $\widehat{\theta}_{2k-1} \geq \cdots \geq \widehat{\theta}_2 > 0$.
4: Evaluate
$$\epsilon_\sigma(\omega_i) = d^{\mathsf{T}} M d \sum_{\substack{j=1 \\ \widehat{\theta}_j > 0}}^{2k-1} |\widehat{S}_k(1,j)|^2 \frac{g_\sigma(\omega_i - \widehat{\theta}_j) - g_\sigma(\omega_i + \widehat{\theta}_j)}{\widehat{\theta}_j}$$
for $i = 1, \ldots, m$.
---

because $d_l \in \mathcal{U}_0$ and $d_r \in \mathcal{U}_\pi$. In the following, we discuss how to preserve this type of structure in a Lanczos procedure.

It was observed in [12] that $H = C_n \Omega$ is self-adjoint with respect to the inner product defined by $\Omega$ in (3), because
$$\langle x, Hy \rangle_\Omega = y^{\mathsf{H}} \Omega C_n \Omega x = \langle Hx, y \rangle_\Omega.$$

We make another observation that $H^2 = (C_n \Omega C_n)\Omega$ is Hermitian and positive definite with respect to the $\Omega$-inner product. Thus there exists a Lanczos procedure associated with $H^2$ that is defined in terms of the $\Omega$-inner product. If we start with the vector $q_1 \in \mathcal{U}_0$, the recurrence relationship among the Lanczos vectors is characterized by the following theorem.

**Theorem 2.** *Let $H = C_n \Omega$ be a definite Bethe–Salpeter Hamiltonian matrix. Suppose that $u_1 \in \mathbb{C}^n$ satisfies $\operatorname{Re}(u_1^{\mathsf{H}} A u_1 + u_1^{\mathsf{H}} B \overline{u}_1) = 1$. Then for $k < n$, applying a $k$-step Lanczos procedure to $H^2$ in the $\Omega$-inner product with the starting vector $[u_1^{\mathsf{H}}, \overline{u}_1^{\mathsf{H}}]^{\mathsf{H}}$ produces*
$$H^2 \begin{bmatrix} U_k \\ \overline{U}_k \end{bmatrix} = \begin{bmatrix} U_k \\ \overline{U}_k \end{bmatrix} T_k + \beta_k \begin{bmatrix} u_{k+1} \\ \overline{u}_{k+1} \end{bmatrix} e_k^{\mathsf{H}}, \tag{21}$$

*where $U_k = [u_1, \ldots, u_k] \in \mathbb{C}^{n \times k}$, $T_k \in \mathbb{R}^{k \times k}$ is as defined in (12). The tridiagonal matrix $T_k$ is positive definite and componentwise nonnegative, and $\beta_k > 0$, if the Lanczos procedure does not break down. The Lanczos vectors satisfy the orthogonality condition*
$$\begin{bmatrix} u_i \\ \overline{u}_i \end{bmatrix}^{\mathsf{H}} \Omega \begin{bmatrix} u_j \\ \overline{u}_j \end{bmatrix} = 2\delta_{ij}, \qquad (1 \leq i, j \leq k+1), \tag{22}$$

where $\delta_{ij}$ is the Kronecker delta notation.

*Proof.* In the generic case (i.e., assuming no breakdown occurs), the Arnoldi procedure using the orthogonality condition (22) with starting vector $q_1 = [u_1^{\mathsf{H}}, \overline{u}_1^{\mathsf{H}}]$ reads

$$H^2 Q_k = Q_k T_k + \beta_k q_{k+1} e_k^{\mathsf{H}},$$

where $T_k$ is an upper Hessenberg matrix with positive subdiagonal entries, and $\beta_k > 0$. Multiplying from the left by $Q_k^{\mathsf{H}} \Omega$, we obtain that

$$2 T_k = Q_k^{\mathsf{H}} \Omega H^2 Q_k = (C_n \Omega Q_k)^{\mathsf{H}} \Omega (C_n \Omega Q_k)$$

is Hermitian positive definite. Consequently the diagonal entries of $T_k$ are real and positive. Hence we conclude that $T_k$ is real symmetric, tridiagonal, positive definite, and componentwise nonnegative. The Arnoldi procedure is in fact a Lanczos procedure.

Let us denote by $\alpha_i$ and $\beta_i$, respectively, the $i$th diagonal and subdiagonal entries of $T_k$, i.e., $T_k$ is of the form (12). Notice that $q_1 \in \mathcal{U}_0$ implies $H^2 q_1 \in \mathcal{U}_0$. From the Lanczos procedure we have

$$q_2 = \frac{1}{\beta_1}(H^2 q_1 - \alpha_1 q_1) \in \mathcal{U}_0,$$

because both $\alpha_1$ and $\beta_1$ are real. By induction, we have

$$q_{i+1} = \frac{1}{\beta_i}(H^2 q_i - \alpha_i q_i - \beta_{i-1} q_{i-1}) \in \mathcal{U}_0$$

for $i = 2, \ldots, k$, as the linear combination on the vectors from $\mathcal{U}_0$ involves only real coefficients. This completes the proof. $\square$

In Section 4.4 we show that this Lanczos procedure reduces to the one given in Section 4.1 for real BSE. The additional factor of two in (22) is introduced to make the two Lanczos procedures identical.

It may appear that the Lanczos procedure associated with $H^2$ only provides one of the two sets of vectors required to construct approximations to the oscillator strength. The following observation shows that the other set of vectors can be easily recovered. Let

$$\begin{bmatrix} V_k \\ \overline{V}_k \end{bmatrix} = \Omega \begin{bmatrix} U_k \\ \overline{U}_k \end{bmatrix}, \qquad (23)$$

or, equivalently,

$$\begin{bmatrix} V_k \\ -\overline{V}_k \end{bmatrix} = H \begin{bmatrix} U_k \\ \overline{U}_k \end{bmatrix}.$$

The $U_k$ and $V_k$ matrices can also be generated together from the following recurrence

$$H \begin{bmatrix} U_k & V_k \\ \overline{U}_k & -\overline{V}_k \end{bmatrix} = \begin{bmatrix} U_k & V_k \\ \overline{U}_k & -\overline{V}_k \end{bmatrix} \begin{bmatrix} 0 & T_k \\ I_k & 0 \end{bmatrix} + \beta_k \begin{bmatrix} u_{k+1} \\ \overline{u}_{k+1} \end{bmatrix} e_{2k}^{\mathsf{H}}. \qquad (24)$$

The orthogonality condition (22) becomes

$$\begin{bmatrix} U_k \\ \overline{U}_k \end{bmatrix}^{\mathsf{H}} \begin{bmatrix} V_k \\ \overline{V}_k \end{bmatrix} = 2I_k. \tag{25}$$

However, this condition is not sufficient for constructing the (oblique) projector associated with the subspace

$$\operatorname{span} \begin{bmatrix} U_k & V_k \\ \overline{U}_k & -\overline{V}_k \end{bmatrix}.$$

We show a stronger result in the following theorem.

**Theorem 3.** *Under the same assumption given in Theorem 2. Let $U_k$ and $V_k$ be as defined in (23) and (24). Then*

$$\begin{bmatrix} V_k & U_k \\ \overline{V}_k & -\overline{U}_k \end{bmatrix}^{\mathsf{H}} \begin{bmatrix} U_k & V_k \\ \overline{U}_k & -\overline{V}_k \end{bmatrix} = 2I_{2k}. \tag{26}$$

*Proof.* Since

$$\begin{bmatrix} V_k & U_k \\ \overline{V}_k & -\overline{U}_k \end{bmatrix}^{\mathsf{H}} \begin{bmatrix} U_k & V_k \\ \overline{U}_k & -\overline{V}_k \end{bmatrix} = \begin{bmatrix} 2I_k & V_k^{\mathsf{H}} V_k - \overline{V_k^{\mathsf{H}} V_k} \\ U_k^{\mathsf{H}} U_k - \overline{U_k^{\mathsf{H}} U_k} & 2I_k \end{bmatrix},$$

it suffices to show that $u_i^{\mathsf{H}} u_j$ and $v_i^{\mathsf{H}} v_j$ are both real for all $i$ and $j$. The proof is based on the fact that

$$\begin{bmatrix} u \\ \pm \overline{u} \end{bmatrix}^{\mathsf{H}} (H^{\mathsf{H}})^{\ell_1} C_n H^{\ell_2} \begin{bmatrix} u \\ \pm \overline{u} \end{bmatrix} = \begin{bmatrix} u \\ \mp \overline{u} \end{bmatrix}^{\mathsf{H}} (C_n \Omega)^{(\ell_1+\ell_2)/2} C_n (\Omega C_n)^{(\ell_1+\ell_2)/2} \begin{bmatrix} u \\ \mp \overline{u} \end{bmatrix} = 0$$

holds for any $u \in \mathbb{C}^n$ and any nonnegative integers $\ell_1$, $\ell_2$ as long as $\ell_1 + \ell_2$ is even.

From (21) it can be verified that $[u_j^{\mathsf{H}}, \overline{u}_j^{\mathsf{H}}]^{\mathsf{H}}$ can be expressed as

$$\begin{bmatrix} u_j \\ \overline{u}_j \end{bmatrix} = p_j(H^2) \begin{bmatrix} u_1 \\ \overline{u}_1 \end{bmatrix},$$

where $p_j(\cdot)$ is a polynomial of degree $j$ with real coefficients. Then we obtain

$$2\mathrm{i} \cdot \mathrm{Im}(u_i^{\mathsf{H}} u_j) = \begin{bmatrix} u_i \\ \overline{u}_i \end{bmatrix}^{\mathsf{H}} C_n \begin{bmatrix} u_j \\ \overline{u}_j \end{bmatrix} = \begin{bmatrix} u_1 \\ \overline{u}_1 \end{bmatrix}^{\mathsf{H}} p_i(H^2)^{\mathsf{H}} C_n p_j(H^2) \begin{bmatrix} u_1 \\ \overline{u}_1 \end{bmatrix} = 0$$

by expanding $p_i(H^2)^{\mathsf{H}} C_n p_j(H^2)$ as the sum of monomials. Similarly, $[v_j^{\mathsf{H}}, -\overline{v}_j^{\mathsf{H}}]^{\mathsf{H}}$ can be expressed as

$$\begin{bmatrix} v_j \\ -\overline{v}_j \end{bmatrix} = H \begin{bmatrix} u_j \\ \overline{u}_j \end{bmatrix} = H p_j(H^2) \begin{bmatrix} u_1 \\ \overline{u}_1 \end{bmatrix} = p_j(H^2) H \begin{bmatrix} u_1 \\ \overline{u}_1 \end{bmatrix} = p_j(H^2) \begin{bmatrix} v_1 \\ -\overline{v}_1 \end{bmatrix},$$

and then

$$2\mathrm{i} \cdot \mathrm{Im}(v_i^{\mathsf{H}} v_j) = \begin{bmatrix} v_i \\ -\overline{v}_i \end{bmatrix}^{\mathsf{H}} C_n \begin{bmatrix} v_j \\ -\overline{v}_j \end{bmatrix} = \begin{bmatrix} v_1 \\ -\overline{v}_1 \end{bmatrix}^{\mathsf{H}} p_i(H^2)^{\mathsf{H}} C_n p_j(H^2) \begin{bmatrix} v_1 \\ -\overline{v}_1 \end{bmatrix} = 0. \qquad \square$$

**Algorithm 4** Lanczos procedure in $\Omega$-inner product for complex full BSE.

**Input:** A definite Bethe–Salpeter Hamiltonian matrix $H \in \mathbb{C}^{2n \times 2n}$; a starting vector $u_1 \in \mathbb{C}^n$ satisfying $\left[u_1^{\mathsf{H}}, \overline{u}_1^{\mathsf{H}}\right]^{\mathsf{H}} \Omega \left[u_1^{\mathsf{H}}, \overline{u}_1^{\mathsf{H}}\right]^{\mathsf{H}} = 2$; the number of Lanczos steps, $k$.

**Output:** $\alpha_1, \ldots, \alpha_k, \beta_1, \ldots, \beta_k \in \mathbb{R}$, and $u_1, \ldots, u_{k+1}, v_1, \ldots, v_{k+1} \in \mathbb{C}^n$ satisfying (23)–(25).

1: $\beta_0 \leftarrow 0, \quad u_0 \leftarrow 0, \quad v_0 \leftarrow 0$.
2: $v_1 = Au_1 + B\overline{u}_1$.
3: **for** $j = 1, \ldots, k$ **do**
4: $\quad x \leftarrow Av_j - B\overline{v}_j - \beta_{j-1}u_{j-1}$.
5: $\quad \alpha_j \leftarrow \operatorname{Re}(v_j^{\mathsf{H}} x)$.
6: $\quad x \leftarrow x - \alpha_j u_j$.
7: $\quad y \leftarrow Ax + B\overline{x}$.
8: $\quad \beta_j \leftarrow \sqrt{\operatorname{Re}(x^{\mathsf{H}} y)}$.
9: $\quad u_{j+1} \leftarrow x/\beta_j, \quad v_{j+1} \leftarrow y/\beta_j$.
10: **end for**

From Theorem 3, we conclude that

$$\frac{1}{2}\begin{bmatrix} U_k & V_k \\ \overline{U}_k & -\overline{V}_k \end{bmatrix} \begin{bmatrix} V_k & U_k \\ \overline{V}_k & -\overline{U}_k \end{bmatrix}^{\mathsf{H}}$$

is the projector we seek, and

$$\begin{bmatrix} 0 & T_k \\ I_k & 0 \end{bmatrix} = \frac{1}{2}\begin{bmatrix} V_k & U_k \\ \overline{V}_k & -\overline{U}_k \end{bmatrix}^{\mathsf{H}} H \begin{bmatrix} U_k & V_k \\ \overline{U}_k & -\overline{V}_k \end{bmatrix}.$$

is indeed a projected form of $H$.

The recurrence given by (24) is more desirable than that given by (21) because it removes the ambiguity introduced by squaring the eigenvalues of the projected matrix

$$\begin{bmatrix} 0 & T_k \\ I_k & 0 \end{bmatrix},$$

which appear in pairs $\pm\theta_i$, where $\theta_i^2$ is the eigenvalue of $T_k$. We regard (24) as a structure preserving Lanczos procedure as the spectrum of the projected matrix is real and symmetric with respect to the origin. Algorithm 4 outlines the structure preserving Lanczos procedure for complex BSE. Similar to Algorithm 2, a careful formulation of short recurrence instead of full orthogonalization is used to largely retain numerical stability.

Finally, we remark that this Lanczos procedure can be extended to have a starting vector from $\mathcal{U}_\phi$. Let $D_\phi = \operatorname{diag}\{I_n, e^{i\phi}I_n\}$. Notice that $D_\phi^{\mathsf{H}} H D_\phi = C_n(D_\phi^{\mathsf{H}} \Omega D_\phi)$ is also a definite BSH matrix. Thus, the Lanczos procedure of $D_\phi^{\mathsf{H}} H D_\phi$,

$$(D_\phi^{\mathsf{H}} H D_\phi) \begin{bmatrix} U_k & V_k \\ \overline{U}_k & -\overline{V}_k \end{bmatrix} = \begin{bmatrix} U_k & V_k \\ \overline{U}_k & -\overline{V}_k \end{bmatrix} \begin{bmatrix} 0 & T_k \\ I_k & 0 \end{bmatrix} + \beta_k \begin{bmatrix} u_{k+1} \\ \overline{u}_{k+1} \end{bmatrix} e_{2k}^{\mathsf{H}},$$

is equivalent to

$$H \begin{bmatrix} U_k & V_k \\ e^{i\phi}\overline{U}_k & -e^{i\phi}\overline{V}_k \end{bmatrix} = \begin{bmatrix} U_k & V_k \\ e^{i\phi}\overline{U}_k & -e^{i\phi}\overline{V}_k \end{bmatrix} \begin{bmatrix} 0 & T_k \\ I_k & 0 \end{bmatrix} + \beta_k \begin{bmatrix} u_{k+1} \\ e^{i\phi}\overline{u}_{k+1} \end{bmatrix} e_{2k}^{\mathsf{H}}.$$

## 4.3 Estimation of the absorption spectrum

In the following we describe how to use the Lanczos procedure defined by (24) to estimate the absorption spectrum. It follows from (24) and the orthogonality condition (26) that

$$\epsilon_\sigma(\omega) = d_r^{\mathsf{H}} g_\sigma(\omega I_{2n} - H) d_l$$

$$= \frac{1}{2} \|d_l\|_\Omega^2 \begin{bmatrix} u_1 \\ -\overline{u}_1 \end{bmatrix}^{\mathsf{H}} g_\sigma(\omega I_{2n} - H) \begin{bmatrix} u_1 \\ \overline{u}_1 \end{bmatrix}$$

$$\approx \frac{1}{4} \|d_l\|_\Omega^2 \begin{bmatrix} u_1 \\ -\overline{u}_1 \end{bmatrix}^{\mathsf{H}} \begin{bmatrix} U_k & V_k \\ \overline{U}_k & -\overline{V}_k \end{bmatrix} g_\sigma\left(\omega I_{2n} - \begin{bmatrix} 0 & T_k \\ I_k & 0 \end{bmatrix}\right) \begin{bmatrix} V_k & U_k \\ \overline{V}_k & -\overline{U}_k \end{bmatrix}^{\mathsf{H}} \begin{bmatrix} u_1 \\ \overline{u}_1 \end{bmatrix}. \quad (27)$$

In the proof of Theorem 3, we showed that $U_k^{\mathsf{H}} u_1$ is real. As a result, we obtain that

$$\begin{bmatrix} u_1 \\ -\overline{u}_1 \end{bmatrix}^{\mathsf{H}} \begin{bmatrix} U_k & V_k \\ \overline{U}_k & -\overline{V}_k \end{bmatrix} = 2 \begin{bmatrix} 0 \\ e_1 \end{bmatrix}^{\mathsf{H}}, \quad \begin{bmatrix} V_k & U_k \\ \overline{V}_k & -\overline{U}_k \end{bmatrix}^{\mathsf{H}} \begin{bmatrix} u_1 \\ \overline{u}_1 \end{bmatrix} = 2 \begin{bmatrix} e_1 \\ 0 \end{bmatrix}.$$

These equations allow us to further simplify the expression given in (27). The simplification removes $U_k$ and $V_k$ in the approximation of $\epsilon_\sigma(\omega)$. Hence these vectors do not need to be explicitly stored. Let $T_k = S_k \Theta_k^2 S_k^{\mathsf{H}}$ be the spectral decomposition of $T_k$, where $\Theta_k = \text{diag}\{\theta_1, \ldots, \theta_k\} \succ 0$. By simple calculation, we obtain

$$f\left(\begin{bmatrix} 0 & T_k \\ I_k & 0 \end{bmatrix}\right) = \begin{bmatrix} S_k & 0 \\ 0 & S_k \end{bmatrix} f\left(\begin{bmatrix} 0 & \Theta_k^2 \\ I_k & 0 \end{bmatrix}\right) \begin{bmatrix} S_k & 0 \\ 0 & S_k \end{bmatrix}^{\mathsf{H}}$$

$$= \frac{1}{2} \begin{bmatrix} S_k & 0 \\ 0 & S_k \end{bmatrix} \begin{bmatrix} \Theta_k & -\Theta_k \\ I_k & I_k \end{bmatrix} \begin{bmatrix} f(\Theta_k) & 0 \\ 0 & f(-\Theta_k) \end{bmatrix} \begin{bmatrix} \Theta_k^{-1} & I_k \\ -\Theta_k^{-1} & I_k \end{bmatrix} \begin{bmatrix} S_k & 0 \\ 0 & S_k \end{bmatrix}^{\mathsf{H}}$$

and

$$\begin{bmatrix} 0 \\ e_1 \end{bmatrix}^{\mathsf{H}} f\left(\begin{bmatrix} 0 & T_k \\ I_k & 0 \end{bmatrix}\right) \begin{bmatrix} e_1 \\ 0 \end{bmatrix} = \frac{1}{2} e_1^{\mathsf{H}} S_k [f(\Theta_k) - f(-\Theta_k)] \Theta_k^{-1} S_k^{\mathsf{H}} e_1 \quad (28)$$

for any smooth function $f(t)$. Substituting $f(t) = f(t;\omega) = g_\sigma(\omega - t)$, we finally arrive at

$$\epsilon_\sigma(\omega) \approx \frac{1}{2} \|d_l\|_\Omega^2 e_1^{\mathsf{H}} S_k [g_\sigma(\omega I_k - \Theta_k) - g_\sigma(\omega I_k + \Theta_k)] \Theta_k^{-1} S_k^{\mathsf{H}} e_1$$

$$= \text{Re}(d^{\mathsf{H}} A d + d^{\mathsf{H}} B \overline{d}) \sum_{j=1}^{k} |S_k(1,j)|^2 \frac{g_\sigma(\omega - \theta_j) - g_\sigma(\omega + \theta_j)}{\theta_j}. \quad (29)$$

**Algorithm 5** The Lanczos algorithm for estimating the optical absorption spectrum for complex full BSE.

---

**Input:** A definite Bethe–Salpeter Hamiltonian matrix $H \in \mathbb{C}^{2n \times 2n}$, an optical transition vector $d \in \mathbb{R}^n$, a broadening factor $\sigma > 0$, the number of Lanczos steps $k$, and a set of frequencies $\{\omega_i\}_{i=1}^m$.

**Output:** The estimated absorption spectrum $\epsilon_\sigma(\omega)$ sampled at $\omega_i$ (for $1 \leq i \leq m$).

1: Perform $k$ Lanczos steps in $\Omega$-inner product with starting vector $[d^{\mathsf{H}}, \overline{d}^{\mathsf{H}}]^{\mathsf{H}}$ using Algorithm 4.
2: Formulate $\widehat{T}_k$ as defined in (14).
3: Compute the spectral decomposition $\widehat{T}_k = \widehat{S}_k \operatorname{diag}\left\{\widehat{\theta}_1^2, \ldots, \widehat{\theta}_{2k-1}^2\right\} \widehat{S}_k^{\mathsf{H}}$, where $\widehat{S}_k^{\mathsf{H}} \widehat{S}_k = I_{2k-1}$ and $\widehat{\theta}_{2k-1} \geq \cdots \geq \widehat{\theta}_2 > 0$.
4: Evaluate

$$\epsilon_\sigma(\omega_i) = \operatorname{Re}\bigl(d^{\mathsf{H}} A d + d^{\mathsf{H}} B \overline{d}\bigr) \sum_{\substack{j=1 \\ \widehat{\theta}_j > 0}}^{2k-1} |\widehat{S}_k(1,j)|^2 \frac{g_\sigma(\omega_i - \widehat{\theta}_j) - g_\sigma(\omega_i + \widehat{\theta}_j)}{\widehat{\theta}_j}$$

for $i = 1, \ldots, m$.

---

Again we have the desired property that $\epsilon_\sigma(\omega) \geq 0$ always holds for $\omega > 0$.

The technique of generalized averaged Gauss quadrature can also be adopted here. Notice that (28) can be interpreted as

$$\begin{bmatrix} 0 \\ e_1 \end{bmatrix}^{\mathsf{H}} f\left(\begin{bmatrix} 0 & T_k \\ I_k & 0 \end{bmatrix}; \omega\right) \begin{bmatrix} e_1 \\ 0 \end{bmatrix} = e_1^{\mathsf{H}} h(T_k; \omega) e_1, \qquad (30)$$

where $h(t; \omega) = t^{-1/2}\bigl[f(t^{1/2}; \omega) - f(-t^{1/2}; \omega)\bigr]$. We expect to obtain a better approximation by replacing $T_k$ in (30) with $\widehat{T}_k$ defined in (14). Certainly, the identity matrix $I_k$ needs to be replaced by $I_{2k-1}$ accordingly. Let $\widehat{T}_k = \widehat{S}_k \widehat{\Theta}_k^2 \widehat{S}_k^{\mathsf{H}}$ be the spectral decomposition of $\widehat{T}_k$, where $\widehat{\Theta}_k = \operatorname{diag}\left\{\widehat{\theta}_1, \ldots, \widehat{\theta}_{2k-1}\right\}$ has at most one nonpositive eigenvalue. The generalized averaged Gauss quadrature produces

$$\epsilon_\sigma(\omega) \approx \operatorname{Re}\bigl(d^{\mathsf{H}} A d + d^{\mathsf{H}} B \overline{d}\bigr) \sum_{\substack{j=1 \\ \widehat{\theta}_j > 0}}^{2k-1} |\widehat{S}_k(1,j)|^2 \frac{g_\sigma(\omega - \widehat{\theta}_j) - g_\sigma(\omega + \widehat{\theta}_j)}{\widehat{\theta}_j}, \qquad (31)$$

which is expected to be better than (29) in general. Algorithm 5 summarizes the Lanczos algorithm with generalized averaged Gauss quadrature for complex full BSE. All of the four desired features listed in the beginning of this section are satisfied.

## 4.4 Connection with other Lanczos procedures

In this subsection, we establish the connection among several variants of the Lanczos procedures. The comparison includes the Lanczos procedures we have discussed in Sections 4.1 and 4.2, as well as that proposed in [33] and [13]. The connection with a variant of the symplectic Lanczos procedure from [38] is also discussed.

**Lanczos procedures for real BSE** In [33, Section 3], a Lanczos procedure that produces

$$\hat{U}_k^\mathsf{T} \hat{V}_k = I_k, \qquad K\hat{U}_k = \hat{V}\hat{T}_k, \qquad M\hat{V}_k = \hat{U}_k \hat{D}_k \tag{32}$$

is studied for real BSE, where $\hat{T}_k$ is symmetric tridiagonal, and $\hat{D}_k \succ 0$ is diagonal. By rescaling $\hat{U}_k$, $\hat{V}_k$ and $\hat{T}_k$ in (32) as

$$U_k = \hat{U}_k \hat{D}_k^{1/2}, \qquad V_k = \hat{V}_k \hat{D}_k^{1/2}, \qquad T_k = \hat{D}_k^{1/2} \hat{T}_k \hat{D}_k^{1/2},$$

we obtain

$$U_k^\mathsf{T} V_k = I_k, \qquad KU_k = VT_k, \qquad MV_k = U_k,$$

which is identical to the Lanczos procedure (18) in the $M$-inner product. As the rescaling is invertible, (32) and (18) are mathematically equivalent. Since there is no need to keep an additional diagonal matrix, (18) is slightly simpler compared to (32).

If both $H$ and the optical transition vector $d$ are real, the Lanczos procedure (24) simplifies to

$$KMU_k = U_k T_k + \beta_k u_{k+1} e_k^\mathsf{H}.$$

The orthogonality condition (25) becomes

$$V_k = MU_k, \qquad V_k^\mathsf{H} U_k = I_k, \tag{33}$$

or simply $U_k^\mathsf{H} M U_k = I_k$. Thus (24) and (18) are identical for real BSE. In the computation of the absorption spectrum for real BSE, (29) and (31) also reduce to (19) and (20), respectively. Therefore, Algorithm 5 can be regarded as a generalization of Algorithm 3 to complex BSE.

**Lanczos procedures for complex BSE** In [12], a Lanczos procedure defined in terms of the $\Omega$-inner product, which produces

$$H\tilde{Q}_k = \tilde{Q}_k \tilde{T}_k + \tilde{\beta}_k \tilde{q}_{k+1} e_k^\mathsf{H}, \qquad \tilde{Q}_{k+1}^\mathsf{H} \Omega \tilde{Q}_{k+1} = I_{k+1}, \tag{34}$$

is presented. However, the projected symmetric tridiagonal matrix $\tilde{T}_k$ does not necessarily have a real spectrum that is symmetric with respect to the origin. Thus (34) is not structured preserving in general. In a subsequent paper [13], it was proposed that a structured starting vector $\tilde{q}_1 \in \mathcal{U}_0$ should be used in (34). With such a structured starting vector, it can be shown that $\tilde{T}_k$ is a real tridiagonal matrix whose diagonal entries are zeros. In addition the nonzero eigenvalues of $\tilde{T}_k$ appear in pairs $\pm\theta$. Hence (34) with $\tilde{q}_1 \in \mathcal{U}_0$ can be regarded as structure preserving. In the following we shall show that this Lanczos procedure is mathematically equivalent to (24).

We have shown that the real symmetric tridiagonal matrix $T_k$ in (21) and (24) is positive definite and componentwise nonnegative. Therefore it admits a Cholesky decomposition $T_k = L_k L_k^{\mathsf{H}}$ where

$$L_k = \mathrm{tridiag}\left\{\begin{matrix} & 0 & \cdots & & 0 & \\ \tilde{\beta}_1 & & \cdots & & & \tilde{\beta}_{2k-1} \\ & \tilde{\beta}_2 & \cdots & & \tilde{\beta}_{2k-2} & \end{matrix}\right\}$$

is a bidiagonal lower triangular matrix, which is also componentwise nonnegative. Multiply $\mathrm{diag}\{I_k, L_k\}^{-\mathsf{H}}$ from the right to (24) yields

$$H\begin{bmatrix} U_k & V_k L_k^{-\mathsf{H}} \\ \overline{U}_k & -\overline{V_k L_k^{-\mathsf{H}}} \end{bmatrix} = \begin{bmatrix} U_k & V_k L_k^{-\mathsf{H}} \\ \overline{U}_k & -\overline{V_k L_k^{-\mathsf{H}}} \end{bmatrix} \begin{bmatrix} 0 & L_k \\ L_k^{\mathsf{H}} & 0 \end{bmatrix} + \beta_k \begin{bmatrix} u_{k+1} \\ \overline{u}_{k+1} \end{bmatrix} \begin{bmatrix} 0 \\ L_k^{-1} e_k \end{bmatrix}^{\mathsf{H}}.$$

Notice that $L_k^{-1} e_k$ is parallel to $e_k$. By setting

$$\tilde{U}_k = \frac{1}{\sqrt{2}} U_k, \qquad \tilde{V}_k = \frac{1}{\sqrt{2}} V_k L_k^{-\mathsf{H}}, \qquad \tilde{\beta}_{2k} = \sqrt{2}\, e_k L_k^{-1} e_k, \qquad (35)$$

we arrive at a Lanczos procedure of the form

$$H\begin{bmatrix} \tilde{U}_k & \tilde{V}_k \\ \overline{\tilde{U}}_k & -\overline{\tilde{V}}_k \end{bmatrix} = \begin{bmatrix} \tilde{U}_k & \tilde{V}_k \\ \overline{\tilde{U}}_k & -\overline{\tilde{V}}_k \end{bmatrix} \begin{bmatrix} 0 & L_k \\ L_k^{\mathsf{H}} & 0 \end{bmatrix} + \tilde{\beta}_{2k} \begin{bmatrix} \tilde{u}_{k+1} \\ \overline{\tilde{u}}_{k+1} \end{bmatrix} e_{2k}^{\mathsf{H}}. \qquad (36)$$

Let

$$\tilde{q}_{2j-1} = \begin{bmatrix} \tilde{u}_{2j-1} \\ \overline{\tilde{u}}_{2j-1} \end{bmatrix}, \qquad \tilde{q}_{2j} = \begin{bmatrix} \tilde{v}_{2j} \\ -\overline{\tilde{v}}_{2j} \end{bmatrix}.$$

Applying the permutation matrix $[e_1, e_{k+1}, e_2, e_{k+2}, \ldots, e_k, e_{2k}]$ from the right to (36) yields

$$H\tilde{Q}_{2k} = \tilde{Q}_{2k} \tilde{T}_{2k} + \tilde{\beta}_k \tilde{q}_{2k+1} e_{2k}^{\mathsf{H}},$$

where

$$\tilde{T}_{2k} = \mathrm{tridiag}\left\{\begin{matrix} & \tilde{\beta}_1 & \tilde{\beta}_2 & \cdots & \tilde{\beta}_{2k-2} & \tilde{\beta}_{2k-1} & \\ 0 & & 0 & \cdots & & 0 & 0 \\ & \tilde{\beta}_1 & \tilde{\beta}_2 & \cdots & \tilde{\beta}_{2k-2} & \tilde{\beta}_{2k-1} & \end{matrix}\right\}.$$

To obtain the orthogonality condition in terms of $\tilde{Q}_k$, we multiply (36) from left by

$$\begin{bmatrix} \tilde{U}_k & \tilde{V}_k \\ \overline{\tilde{U}}_k & -\overline{\tilde{V}}_k \end{bmatrix}^{\mathsf{H}} C_n.$$

Using (26) and simple algebraic manipulation, we obtain

$$\begin{bmatrix} \tilde{U}_k & \tilde{V}_k \\ \overline{\tilde{U}}_k & -\overline{\tilde{V}}_k \end{bmatrix}^{\mathsf{H}} \Omega \begin{bmatrix} \tilde{U}_k & \tilde{V}_k \\ \overline{\tilde{U}}_k & -\overline{\tilde{V}}_k \end{bmatrix} = I_{2k}.$$

Thus we have derived (34) from (24), assuming the number of Lanczos steps in (34) is even. As the transformation (35) is invertible, the two Lanczos procedures are mathematically equivalent.

The Lanczos procedure (34) can be used to approximate the absorption spectrum as follows:

$$\begin{aligned} \epsilon_\sigma(\omega) &= d_r^{\mathsf{H}} g_\sigma(\omega I_{2n} - H) d_l \\ &\approx d_r^{\mathsf{H}} \tilde{Q}_{2k} g_\sigma(\omega I_{2n} - \tilde{T}_{2k}) \tilde{Q}_{2k}^{\mathsf{H}} \Omega d_l \\ &= \frac{1}{2} \|d_l\|_\Omega^2 e_1^{\mathsf{H}} g_\sigma(\omega I_{2n} - \tilde{T}_{2k}) \tilde{T}_{2k}^{-1} e_1. \end{aligned} \tag{37}$$

The derivation of the last step requires similar effort compared to the proof of Theorem 3. The expression (37) is also mathematically equivalent to (29). The main difference between them is that the spectral decomposition of $\tilde{T}_{2k}$ instead of that of $T_k$ is needed. However, we remark that there exist subtle differences when the technique generalized averaged Gauss quadrature is adopted. A direct application of generalized averaged Gauss quadrature replaces $\tilde{T}_{2k}$ by a $(4k-1) \times (4k-1)$ tridiagonal matrix

$$\widehat{\tilde{T}}_{2k} = \operatorname{tridiag} \left\{ \begin{array}{ccccccccc} \tilde{\beta}_1 & \cdots & & \tilde{\beta}_{2k-1} & \tilde{\beta}_{2k} & \tilde{\beta}_{2k-2} & \cdots & & \tilde{\beta}_1 \\ 0 & \cdots & \cdots & & 0 & 0 & \cdots & \cdots & 0 \\ \tilde{\beta}_1 & \cdots & & \tilde{\beta}_{2k-1} & \tilde{\beta}_{2k} & \tilde{\beta}_{2k-2} & \cdots & & \tilde{\beta}_1 \end{array} \right\}.$$

The positive eigenvalues of $\widehat{\tilde{T}}_{2k}$ are not quite the same as the those of $\widehat{T}_k^{1/2}$, although the number of positive Gauss nodes in the generalized averaged Gauss quadrature is $2k-1$ for both case. We shall see from the numerical experiments the generalized averaged Gauss quadrature based on $\widehat{\tilde{T}}_{2k}$ is in general slightly worse than that based on $\widehat{T}_k$ in terms of accuracy.

We remark that in the discussion above we always assume that an even number of Lanczos steps is performed in (34). In fact, for an odd number of Lanczos steps, $\tilde{T}_{2k+1}$ always has a zero eigenvalue. In the view of Gauss quadrature for estimating the absorption spectrum, such a zero eigenvalue is not a very useful Gauss quadrature node because $\epsilon_\sigma(0) = 0$ is known trivially. Therefore, an even number of Lanczos steps should be performed when computing (34). Similarly, the zero eigenvalue of $\widehat{\tilde{T}}_{2k}$ is not very helpful. Thus we only consider the $2k-1$ positive eigenvalues of $\widehat{\tilde{T}}_{2k}$ to be useful in the generalized averaged Gauss quadrature.

**Connection with symplectic Lanczos procedure** We have shown that our new Lanczos procedure (24) is essentially equivalent to the one proposed in [13], and both are equivalent to (18) and the one in [33] when applied to real BSE. There exists other equivalent formulations. We present these formulations in this section, and exploit more properties of the Lanczos procedure.

Let

$$\tilde{X}_k = \frac{U_k + V_k}{2}, \qquad \tilde{Y}_k = \frac{\overline{U}_k - \overline{V}_k}{2}, \qquad \tilde{A}_k = \frac{I_k + T_k}{2}, \qquad \tilde{B}_k = \frac{I_k - T_k}{2}.$$

Then we reformulate (24) as

$$H \begin{bmatrix} \tilde{X}_k & \overline{\tilde{Y}}_k \\ \tilde{Y}_k & \overline{\tilde{X}}_k \end{bmatrix} = \begin{bmatrix} \tilde{X}_k & \overline{\tilde{Y}}_k \\ \tilde{Y}_k & \overline{\tilde{X}}_k \end{bmatrix} \begin{bmatrix} \tilde{A}_k & \tilde{B}_k \\ -\tilde{B}_k & -\tilde{A}_k \end{bmatrix} + \frac{1}{2}\beta_k \begin{bmatrix} \tilde{x}_{k+1} & \overline{\tilde{y}}_{k+1} \\ \tilde{y}_{k+1} & \overline{\tilde{x}}_{k+1} \end{bmatrix} \begin{bmatrix} 0 & e_k^{\mathsf{H}} \\ 0 & e_k^{\mathsf{H}} \end{bmatrix}. \qquad (38)$$

The orthogonality condition (26) becomes

$$\left( C_n \begin{bmatrix} \tilde{X}_k & \overline{\tilde{Y}}_k \\ \tilde{Y}_k & \overline{\tilde{X}}_k \end{bmatrix} C_k \right)^{\mathsf{H}} \begin{bmatrix} \tilde{X}_k & \overline{\tilde{Y}}_k \\ \tilde{Y}_k & \overline{\tilde{X}}_k \end{bmatrix} = \begin{bmatrix} \tilde{X}_k & -\overline{\tilde{Y}}_k \\ -\tilde{Y}_k & \overline{\tilde{X}}_k \end{bmatrix}^{\mathsf{H}} \begin{bmatrix} \tilde{X}_k & \overline{\tilde{Y}}_k \\ \tilde{Y}_k & \overline{\tilde{X}}_k \end{bmatrix} = I_{2k}. \qquad (39)$$

Although (24) is derived from (21), which uses the $\Omega$-inner product, the equivalent formulation (38) is a Lanczos procedure in the $C$-inner product. As a result, the projected matrix is a $2k \times 2k$ BSH matrix. As we have discussed in Section 2, the eigenvectors of $H$ are orthogonal in both the $\Omega$-inner product and the $C$-inner product. This suggests that our Lanczos procedure largely preserves properties of $H$. As a byproduct of this observation, we obtain the Cauchy interlacing property as stated in Theorem 4, which provides an estimate on the location of quadrature nodes in the Gauss quadrature. This can be viewed as a generalization of [33, Lemma 3.5]. A proof of Theorem 4 can be found in [29].

**Theorem 4.** *Let $T_k$ be defined as in (21) and suppose that the eigenvalues of $H$ and $T_k$ are $\pm\lambda_1$, $\pm\lambda_2$, ..., $\pm\lambda_n$, and $\theta_1^2$, $\theta_2^2$, ..., $\theta_k^2$, respectively, with $0 < \lambda_1 \leq \lambda_2 \leq \cdots \leq \lambda_n$, $0 < \theta_1 \leq \theta_2 \leq \cdots \leq \theta_k$. Then, under the assumption given in Theorem 2, we have*

$$\lambda_i \leq \theta_i \leq \lambda_{n-k+i}, \qquad (1 \leq i \leq k).$$

It has been shown in [27] that the matrix

$$\mathrm{i} Q_n^{\mathsf{H}} H Q_n = J_n \begin{bmatrix} \mathrm{Re}(A+B) & \mathrm{Im}(A-B) \\ -\mathrm{Im}(A+B) & \mathrm{Re}(A-B) \end{bmatrix} =: J_n \tilde{M}$$

is a real Hamiltonian matrix with $\tilde{M} \succ 0$, where

$$J_n = \begin{bmatrix} 0 & I_n \\ -I_n & 0 \end{bmatrix}, \qquad Q_n = \frac{1}{\sqrt{2}} \begin{bmatrix} I_n & -\mathrm{i}I_n \\ I_n & \mathrm{i}I_n \end{bmatrix}.$$

Using the same unitary transformation, the Lanczos factorization (38) becomes

$$J_n\tilde{M} \begin{bmatrix} \text{Re}(U_k) & \text{Im}(V_k) \\ -\text{Im}(U_k) & \text{Re}(V_k) \end{bmatrix} = \begin{bmatrix} \text{Re}(U_k) & \text{Im}(V_k) \\ -\text{Im}(U_k) & \text{Re}(V_k) \end{bmatrix} \begin{bmatrix} 0 & T_k \\ -I_k & 0 \end{bmatrix} + [\text{rank 1}], \qquad (40)$$

and the orthogonality condition (39) becomes

$$\begin{bmatrix} \text{Re}(U_k) & \text{Im}(V_k) \\ -\text{Im}(U_k) & \text{Re}(V_k) \end{bmatrix}^{\mathsf{H}} J_n \begin{bmatrix} \text{Re}(U_k) & \text{Im}(V_k) \\ -\text{Im}(U_k) & \text{Re}(V_k) \end{bmatrix} = J_k. \qquad (41)$$

Since (41) indicates that the Lanczos vectors associated with $J_n\tilde{M}$ are symplectic, the Lanczos factorization (40) yields in fact a symplectic Lanczos procedure (see, e.g., [4, 6, 38]) for the real Hamiltonian matrix $J_n\tilde{M}$. Such a variant of symplectic Lanczos procedures has been discussed in [38]. Therefore, (24) can also be interpreted as a variant of symplectic Lanczos procedure. Finally we remark that for the purpose of computing the absorption spectrum, the starting vector in (40) should be chosen parallel to $\left[\text{Re}(d)^{\mathsf{H}}, -\text{Im}(d)^{\mathsf{H}}\right]^{\mathsf{H}}$ if (40) is adopted.

### 4.5 Structure preserving Lanczos algorithm with paired starting vectors

Besides several equivalent structure preserving Lanczos procedures, there are also other structure preserving Lanczos procedures with paired starting vectors. Actually, when $\|u_1\|_2 \neq \|v_1\|_2$, Lanczos procedures of the form

$$H \begin{bmatrix} U_k & \overline{V}_k \\ V_k & \overline{U}_k \end{bmatrix} = \begin{bmatrix} U_k & \overline{V}_k \\ V_k & \overline{U}_k \end{bmatrix} \begin{bmatrix} A_k & B_k \\ -\overline{B}_k & -\overline{A}_k \end{bmatrix} + \begin{bmatrix} u_{k+1} & \overline{v}_{k+1} \\ v_{k+1} & \overline{u}_{k+1} \end{bmatrix} \begin{bmatrix} \beta_k e_k^{\mathsf{H}} & 0 \\ 0 & -\overline{\beta}_k e_k^{\mathsf{H}} \end{bmatrix} \qquad (42)$$

can be constructed, where $A_k$ and $B_k$ are tridiagonal, and the orthogonality condition on the Lanczos vectors is either[2]

$$\begin{bmatrix} u_i & \overline{v}_i \\ v_i & \overline{u}_i \end{bmatrix}^{\mathsf{H}} \Omega \begin{bmatrix} u_j & \overline{v}_j \\ v_j & \overline{u}_j \end{bmatrix} = \delta_{ij} I_2 \qquad (43)$$

or

$$C_2 \begin{bmatrix} u_i & \overline{v}_i \\ v_i & \overline{u}_i \end{bmatrix}^{\mathsf{H}} C_n \begin{bmatrix} u_j & \overline{v}_j \\ v_j & \overline{u}_j \end{bmatrix} = \begin{bmatrix} u_i & -\overline{v}_i \\ -v_i & \overline{u}_i \end{bmatrix}^{\mathsf{H}} \begin{bmatrix} u_j & \overline{v}_j \\ v_j & \overline{u}_j \end{bmatrix} = \delta_{ij} I_2. \qquad (44)$$

In fact, from the discussion in the previous subsection, we also see that the condition (44) for BSH matrices is equivalent to the symplecticity condition for real Hamiltonian matrices. For both orthogonality conditions, the eigenvalues of the projected matrix

$$H_k = \begin{bmatrix} A_k & B_k \\ -\overline{B}_k & -\overline{A}_k \end{bmatrix}$$

---

[2] If (43) is used, orthogonalization within each two dimensional subspace is required.

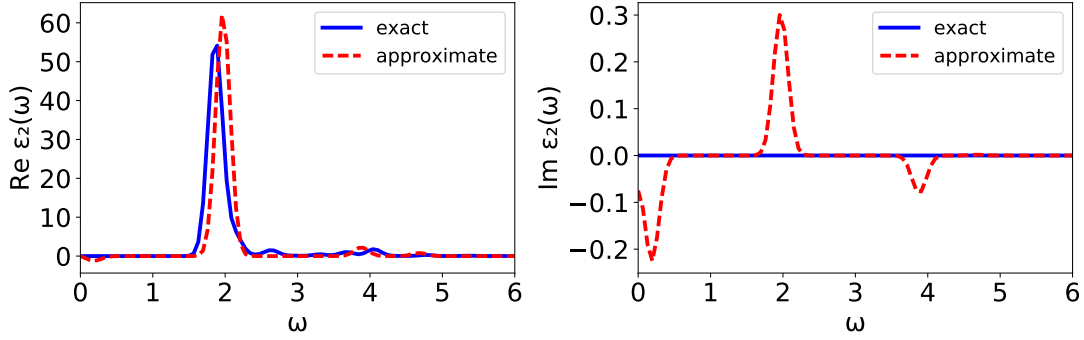

Figure 2: The real (left) and imaginary (right) parts of the absorption spectrum approximated using (45) on a small example provided in (46) with $k = 4$ and $\sigma = 0.1$. The real part is not nonnegative, and the imaginary part is nonzero.

are real and occur in pairs $\pm \theta$.

When estimating the absorption spectrum using (42), we use $u_1 = d$, $v_1 = 0$ as the starting vectors because $d_l$ does not satisfy the condition $\|u_1\|_2 \neq \|v_1\|_2$. We point out that the absorption spectrum computed based on (42) is *not* necessarily real and nonnegative for $\omega > 0$. For instance, if the orthogonality condition (43) is used, the absorption spectrum approximated by

$$\epsilon_\sigma(\omega) \approx d_r^{\mathsf{H}} \begin{bmatrix} U_k & \overline{V}_k \\ V_k & \overline{U}_k \end{bmatrix} g_\sigma(\omega I_k - H_k) \begin{bmatrix} U_k & \overline{V}_k \\ V_k & \overline{U}_k \end{bmatrix}^{\mathsf{H}} \Omega d_l \qquad (45)$$

is not guaranteed to be real, although the Ritz values are always real and appear in pairs. We illustrate this by a small artificial example with $n = 16$:

$$A(i,j) = \begin{cases} 4, & \text{if } i = j, \\ 1, & \text{if } |i-j| = 1, \\ 0, & \text{otherwise,} \end{cases} \qquad B(i,j) = \begin{cases} \mathrm{i}^{i-1} & \text{if } i = j, \\ 0, & \text{otherwise,} \end{cases} \qquad d(i) = (-1)^{i-1}. \quad (46)$$

As shown in Figure 2, after four steps of Lanczos procedure the absorption spectrum computed by (45) with $\sigma = 0.1$ is clearly not real, and the real part is not nonnegative. Hence, we do not consider using (42) in the $\Omega$-inner product and (45) for computing the absorption spectrum.

If the orthogonality condition (44) is adopted, also using $u_1 = d$, $v_1 = 0$ as the starting vectors, the projected matrix $H_k$ is a definite BSH matrix. Let the spectral decomposition of $H_k$ be

$$H_k = \begin{bmatrix} S_1 & \overline{S}_2 \\ S_2 & \overline{S}_1 \end{bmatrix} \begin{bmatrix} \Theta & 0 \\ 0 & -\Theta \end{bmatrix} \begin{bmatrix} S_1 & -\overline{S}_2 \\ -S_2 & \overline{S}_1 \end{bmatrix}^{\mathsf{H}},$$

where $\Theta = \text{diag}\{\theta_1, \ldots, \theta_k\} \succ 0$. It can be shown that (42) and (44) lead to a structure preserving algorithm because

$$\epsilon_\sigma(\omega) \approx \|d\|_2^2 (e_1 - e_{k+1})^{\mathsf{H}} g_\sigma(\omega I - H_k)(e_1 + e_{k+1})$$
$$= \|d\|_2^2 \sum_{j=1}^{k} |S_1(1,j) - S_2(1,j)|^2 [g_\sigma(\omega - \theta_j) - g_\sigma(\omega + \theta_j)]. \tag{47}$$

This formulation possesses the second and third features listed in the beginning of this section. However, theoretically Lanczos procedures in the $C$-inner product may sometimes break down due to $C$-neutral vectors.[3] Such a breakdown is not a lucky breakdown. It is also not very clear how to incorporate the technique of generalized averaged Gauss quadrature in (47).

We remark that in general (47) is not as good as (29) even if generalized averaged Gauss quadrature is not used. Our numerical experiments suggest that (47) typically requires about twice as many as Lanczos steps to achieve the same accuracy level compared to (29).[4] A brief explanation is that for the same number of Lanczos steps $k$, $H$ has been raised to the power $H^{2k}$ in (24), while $H$ has only been raised to the power $H^k$ in (42). A higher polynomial degree potentially provides better approximation quality.

## 5  Computational examples

In this section we present several examples to demonstrate the accuracy and efficiency of the Lanczos algorithm for computing the optical absorption spectrum. We implemented the Lanczos algorithms in software packages BSEPACK [28] and BerkeleyGW [9]. All tests were performed on the Linux cluster Edison at the National Energy Research Scientific Computing Center (NERSC).[5] Each computational node on Edison consists of 64 GB DDR3 1866 MHz memory and two sockets, with a 12-core Intel "Ivy Bridge" processor at 2.4 GHz on each socket. The computational nodes are connected by a Cray Aries network with Dragonfly topology, with 23.7 TB/s global bandwidth. Our tests make use of 10 computational nodes and 24 MPI processes per node. The Fortran 90 implementation of algorithms is compiled by the Intel Fortran compiler, and linked with the Cray LibSci and Cray MPI libraries. No multithreading feature is utilized.

For our calculations, we use a benchmark system consisting of a single-wall $(8,0)$ carbon nanotube with 32 atoms, 128 electrons, and 64 Kohn–Sham spin-degenerate bands in the unit cell. As depicted in Figure 3, this system is periodic along the "c" axis, but confined along the other directions labeled by the axes "a" and "b", which makes this an interesting benchmark system. In particular, as we will discuss, the TDA may or may not be a good approximation depending on the direction of optical excitation in this particular system.

---

[3] A $C$-neutral vector is a vector $v \in \mathbb{C}^{2n}$ which satisfies $v^{\mathsf{H}} C_n v = 0$.
[4] A similar behavior has been observed in [36] when solving the linear response eigenvalue problem.
[5] http://www.nersc.gov/users/computational-systems/edison/

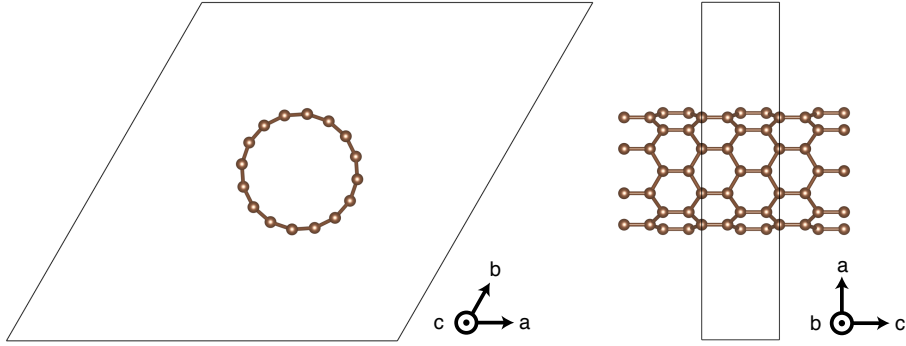

Figure 3: Single-wall (8,0) carbon nanotube benchmark system. The black region represents the unit cell of the system, which is periodic along the "c" axis.

In general, crystal states can be written in a Bloch form as $\Psi_{nk}(r) = e^{ik\cdot r} u_{nk}(r)$, where $n$ is a band index, $k$ is a $k$-point, and $u_{nk}(r)$ is a cell-periodic complex-valued function. Because "c" is the only periodic direction, we only need to sample $k$-points along that axis. When solving the BSE, we include $n_v = 10$ valence states, $n_c = 12$ conduction states, and $n_k = 256$ $k$-points, so that $n = n_v n_c n_k = 30{,}720$, and we picked $g_\sigma$ as a Gaussian function with $\sigma = 100$ meV. The matrices $A$ and $B$ are both dense. We did not perform a systematic convergence test with respect to the number of conducting bands. However, the use of $n_c = 12$ conducting bands already produces main physical features in the absorption spectrum also observed when a larger number of conduction bands are used.

**Full BSE vs. TDA**  In our first experiment, we calculate the absorption spectrum for two different directions for the light polarizations using both full BSE and TDA solvers. To exclude other sources of errors, we fully diagonalize the matrices $H$ and $A$, with dimensions 61,440 and 30,720, respectively. We can see from Figure 4 that even within the same system, the TDA can either be a valid approximation or give a qualitatively wrong absorption spectrum depending on the polarization direction of light. When the polarization of the optical excitation is along the "c" axis, which is a direction along which the system is periodic, the TDA is a good approximation for low-energy optical spectrum. However, if the light polarization is along any confined direction spanned by the "a" and "b" axes, a large difference between the two spectra can be observed. This can be understood from a large exciton–plasmon hybridization which couples to light polarized along the confined direction, and which can not be well-described within the TDA [12].

Thus this example confirms the necessity of developing full BSE solvers for absorption spectrum calculation. In the subsequent tests the light polarizations is chosen to be perpendicular to the tube so that using a full BSE solver is necessary.

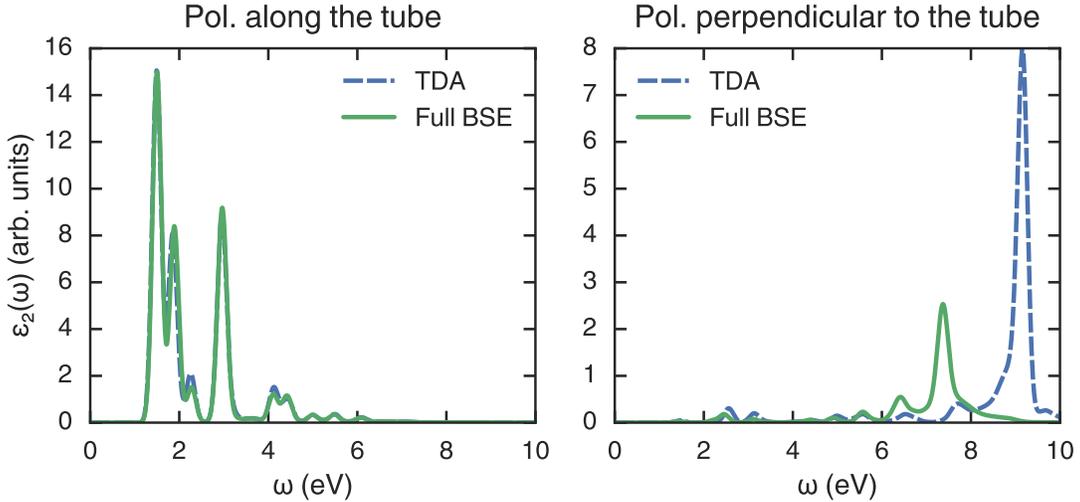

Figure 4: The absorption spectrum for a single-wall $(8,0)$ carbon nanotube for two different directions for the light polarizations. The label "Full BSE" refers to spectrum obtained from solving the full BSH in (8), and "TDA" refers to the spectrum obtained within the Tamm–Dancoff approximation.

**Effectiveness of the Lanczos algorithm** In Figure 5 we plot the approximate absorption spectra obtained by running 32 steps of different variants of the Lanczos algorithm. We use the result obtained from full diagonalization as the "exact" solution to measure the accuracy of these Lanczos algorithms. The paired Lanczos algorithm described in Section 4.5 (Figure 5(a), abbreviated as PL) is clearly worse than the one with a single structured starting vector (Figure 5(b), abbreviated as SVL). The Lanczos algorithm proposed in [13] (abbreviated as GMG) is equivalent to ours with a single structured starting vector and is hence omitted here. The technique of generalized averaged Gauss quadrature (abbreviated as GAGQ) clearly improves the accuracy of the Lanczos algorithm. With such a small number of Lanczos steps, our Lanczos algorithm with generalized averaged Gauss quadrature (i.e., Algorithm 5) already produces very satisfactory result. Though not very clear from this figure, our Algorithm 5 (Figure 5(d)) is slightly better than the one proposed in [13] with generalized averaged Gauss quadrature (Figure 5(c)) in this example.

To measure the accuracy of approximate absorption spectrum, we introduce the concept of *angle* between two functions as follows. Let $\xi(\omega)$ and $\zeta(\omega)$ be sufficiently smooth functions of $\omega$ over an interval $I$. Then the angle between $\xi(\omega)$ and $\zeta(\omega)$ is defined as

$$\angle(\xi, \eta) = \arccos \frac{\langle \xi, \zeta \rangle}{\sqrt{\langle \xi, \xi \rangle \langle \zeta, \zeta \rangle}}, \tag{48}$$

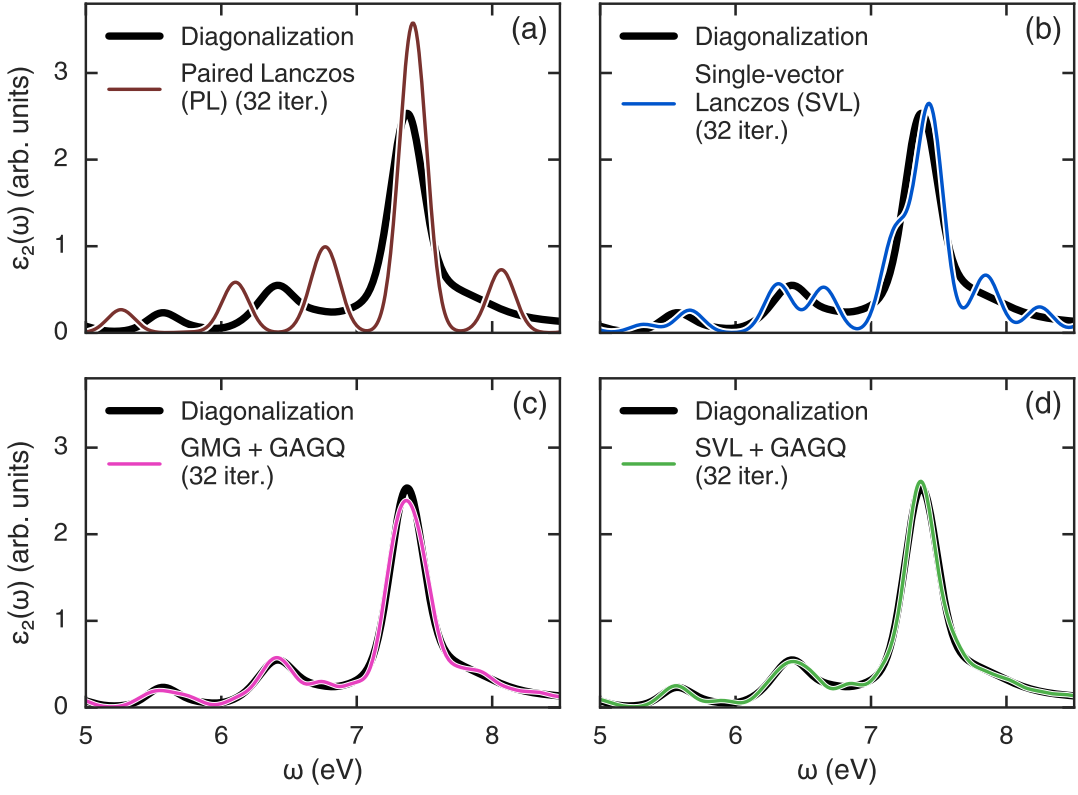

Figure 5: Comparison of the absorption spectra obtained from different variants of the Lanczos algorithm with the spectrum obtained from the full diagonalization of the BSH matrix.

where
$$\langle \xi, \zeta \rangle = \int_I \xi(\omega)\overline{\zeta(\omega)}\,\mathrm{d}\omega$$

is the usual $L^2$-inner product. The angle $\angle(\xi, \eta)$ is in fact the principal angle (also known as canonical angle) between two subspaces, $\mathrm{span}\{\xi(\omega)\}$ and $\mathrm{span}\{\eta(\omega)\}$, of $L^2(I)$. A small angle between two functions implies similar shapes of their curves. This allows us to measure the error of the approximate absorption spectrum compared to the "exact" one in terms of the angle between them. This measure is similar to the cross-correlation measure between two curves. In Figure 6 we plot the errors of different variants of the Lanczos algorithm, with the integrals in (48) approximated by rectangular rules using the sampling points of $\epsilon_2(\omega)$'s. It confirms our observation from Figure 5, not only for a single snapshot after 32 Lanczos steps, but also consistently throughout the whole iterative procedure. The difference between the two different variants with generalized averaged Gauss quadrature becomes more clear in Figure 6. Overall Algorithm 5 is better than the variant from [13] combined with generalized averaged Gauss quadrature. We remark that there are about 10% cases in this example involving nonpositive definite $\widehat{T}_k$ in Algorithm 5. Figure 6 shows that dropping the nonpositive eigenvalue of $\widehat{T}_k$ does not harm the accuracy.

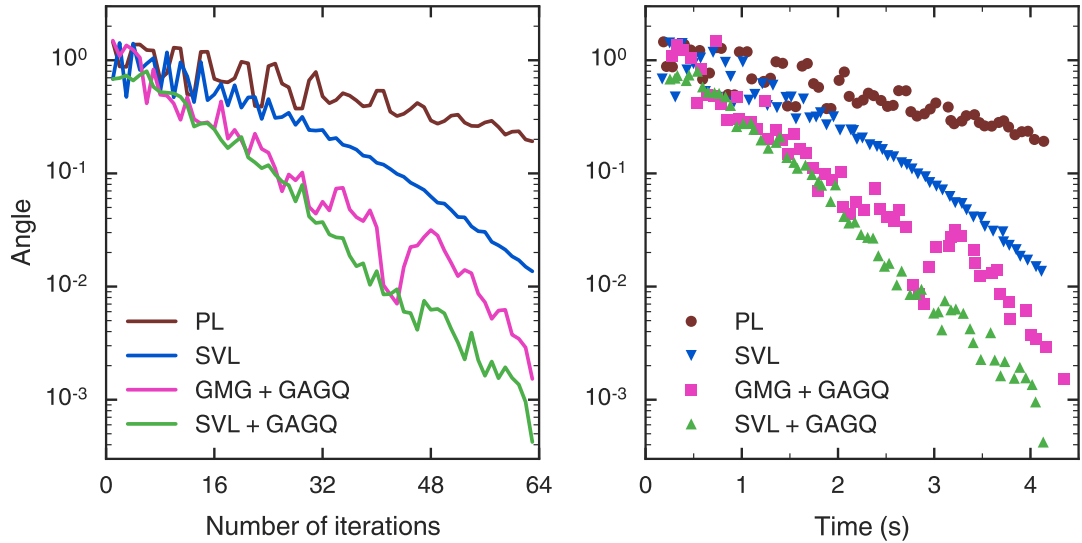

Figure 6: The convergence history of different variants of the Lanczos algorithm. The error is measured by the angle (48) between the approximate absorption spectrum and the one obtained from full diagonalization. The algorithms used here are the same as those in Figure 5.

It takes 62 iterations and 4.1 seconds for Algorithm 5 to achieve the accuracy level $10^{-3}$ (in terms of angles), which is more than sufficient for practical use. This is over 500 times faster compared to full diagonalization (2125.8 seconds). If the multiplications of $A$ and $B$ with vectors can be implemented more efficiently by further exploiting the structures of $A$ and $B$ (see, e.g., [3]), the improvement is expected to be more significant.

In our test, the number of Lanczos steps is always prescribed by the user. We remark that it is possible to instead specify the desired accuracy in the input, and automatically determine the required number of Lanczos steps in the calculation. One strategy proposed in [19] is to estimate the error using the difference between the results obtained with and without generalized averaged Gauss quadrature. However, since this strategy relies on the result without generalized averaged Gauss quadrature, which is in a relatively low accuracy as we have shown in Figure 6, the estimate is in general too pessimistic. A better strategy is to use the difference between two consecutive iterations (i.e., $(k-1)$th and $k$th steps) instead in the stopping criterion.

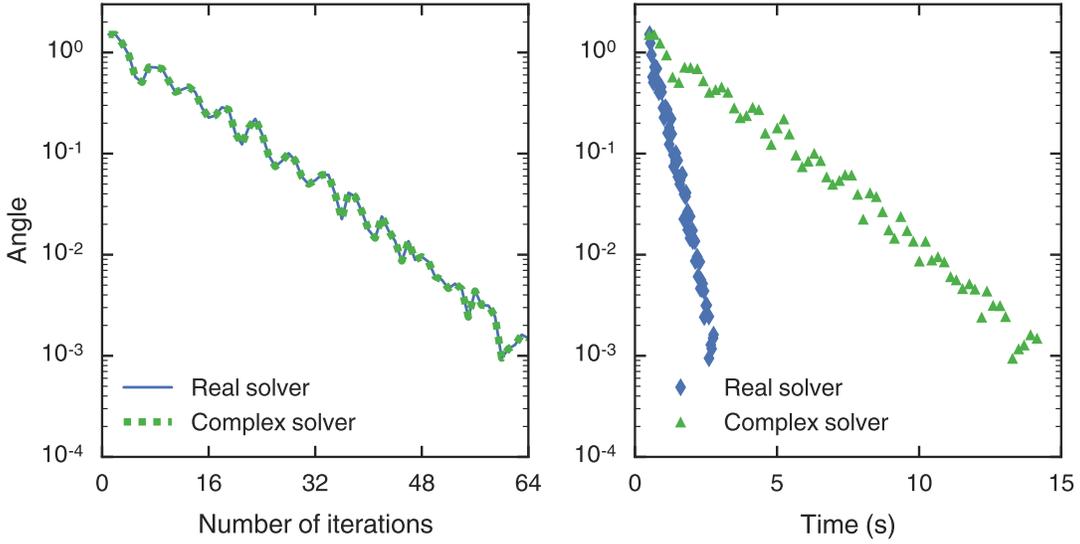

Figure 7: The convergence history of Algorithms 3 and 5 for bulk silicon. The error is measured by the angle (48) between the approximate absorption spectrum and the one obtained from full diagonalization.

**Systems with real-space inversion symmetry** Our last example uses another system which has real-space inversion symmetry. When a system has real-space and time-reversal symmetry, the wave functions in reciprocal space, and thus the BSH matrix, can be written as real numbers [9]. We use bulk silicon for this benchmark, with $n_v = 4$, $n_c = 6$, and $n_k = 1{,}000$, so that the dimension of the $A$ and $B$ blocks of BSH is $n = 24{,}000$. We also use a Gaussian broadening in this system, but with $\sigma = 150$ meV. Since the BSH matrix is real, both Algorithms 3 and 5 are applicable, and are identical as discussed in Section 4. Our experiment confirms the theoretical prediction. In Figure 7 we plot the convergence history (in terms of angles) of the two algorithms. The curves of convergence history are indeed on the top of each other in the left plot. As for the execution time, the real solver is faster than the complex one, due to some additional operations involving the imaginary parts in the complex solver when applying to real matrices.

## 6 Concluding remarks

In this paper we presented and analyzed a simple structure preserving Lanczos procedure for definite Bethe–Salpeter Hamiltonian matrices and combined it with the recently developed technique of generalized averaged Gauss quadrature to estimate the optical absorption spectrum. The analyzed Lanczos procedure possesses several attractive features, such as being free of serious breakdown, and preserving nonnegativity of the absorption spectrum.

The use of alternative inner products based on the orthogonalities of the eigenvectors plays a key role in preserving the structure. By some theoretical analysis we established the equivalence between the Lanczos procedure presented in this paper with several existing Lanczos procedures in the literature, including the ones in [13, 33] for random phase approximation, and one variant of symplectic Lanczos procedure in [38]. Numerical experiments demonstrate that the Lanczos algorithm can provide accurate approximation of the absorption spectrum with a relatively small number of Lanczos steps. In addition, the technique of generalized averaged Gauss quadrature largely improves the accuracy of the Lanczos algorithm. When this technique is applied, our Lanczos algorithm is more efficient and more accurate compared to other variants.

In this work the blocks $A$ and $B$ in the BSH matrix $H$ are formed as dense matrices. However, the Lanczos algorithm does not require these matrices to be explicitly formed. An implicit representation that allows to perform matrix–vector multiplication suffices. Efficient ways of constructing and applying the BSH matrix have been described in [16, 17, 22, 23, 24]. The development of other approximation and compression strategies is planned as future work.

## Acknowledgments

We thank Peter Benner for pointing out the equivalence between (24) and the symplectic Lanczos procedure discussed in [38] for Hamiltonian matrices, and for kindly sharing with us several unpublished notes of his early work in this direction. We also thank Jiri Brabec and Niranjan Govind for helpful discussions.